\documentclass[reqno,11pt]{amsart}

\usepackage[top=2.0cm,bottom=2.0cm,left=3cm,right=3cm]{geometry}
\usepackage{amsthm,amsmath,amssymb,dsfont}
\usepackage{mathrsfs,amsfonts,functan,extarrows,mathtools}
\usepackage[colorlinks]{hyperref}

\usepackage{marginnote}
\usepackage{xcolor}
\usepackage{stmaryrd}
\usepackage{esint}
\usepackage{graphicx}
\usepackage{bbm}

\usepackage{tikz}

\newtheorem{theorem}{Theorem}[section]
\newtheorem{proposition}{Proposition}[section]

\newtheorem{remark}{Remark}[section]
\newtheorem{lemma}{Lemma}[section]

\numberwithin{equation}{section}
\allowdisplaybreaks

\def\p{\partial}
\def\d{\mathrm{d}}
\def\no{\nonumber}
\def\R{\mathbb{R}}
\def\eps{\varepsilon}
\def\div{\mathrm{div}}
\def\u{\mathbf{u}}

\def\l{\langle}
\def\r{\rangle}
\def\exp{\mathrm{exp}}
\def\M{\mathfrak{M}}
\def\A{\mathcal{A}}
\def\B{\mathcal{B}}

\def\u{\mathfrak{u}}
\def\L{\mathcal{L}}

\def\P{\mathcal{P}}
\def\I{\mathcal{I}}

\newcounter{wronumber}

\begin{document}
\title[Compressible Euler limit from Boltzmann equation]
			{Compressible Euler limit from Boltzmann equation with complete diffusive boundary condition in half-space}

\author[Ning Jiang]{Ning Jiang}
\address[Ning Jiang]{\newline School of Mathematics and Statistics, Wuhan University, Wuhan, 430072, P. R. China}
\email{njiang@whu.edu.cn}

\author[Yi-Long Luo]{Yi-Long Luo${}^*$}
\address[Yi-Long Luo]
{\newline School of Mathematics, South China University of Technology, Guangzhou, 510641, P. R. China}
\email{luoylmath@scut.edu.cn}

\author[Shaojun Tang]{Shaojun Tang}
\address[Shaojun Tang]
{\newline Department of Mathematics, Wuhan University of Technology, Wuhan, 430070, P. R. China }
\email{shaojun.tang@whut.edu.cn}
\thanks{${}^*$ Corresponding author \quad \today}

\maketitle

\begin{abstract}
   In this paper, we prove the compressible Euler limit from the Boltzmann equation with hard sphere collisional kernel and complete diffusive boundary condition in half-space by employing the Hilbert expansion which includes interior and Knudsen layers. This rigorously justifies the corresponding formal analysis in Sone's book \cite{Sone-2007-Book} in the context of short time smooth solutions, and also generalizes the classic Caflisch's result \cite{Caflish-1980-CPAM} to initial-boundary problem case. \\


\end{abstract}





\section{Introduction}

\subsection{Boltzmann equation with complete diffusive boundary condition}

In this paper, the following Boltzmann equation with Euler scaling over $(t, x, v) \in \R_+ \times \R^3_+ \times \R^3$ is considered:
\begin{equation}\label{BE}
\left\{
\begin{array}{l}
\p_t F_\eps + v \cdot \nabla_x F_\eps = \frac{1}{\eps} B (F_\eps, F_\eps) \quad \qquad {\text{on}}\  \R_+ \times \R^3_+ \times \R^3 \,,\\[1.5mm]
F_\eps (0, x, v) = F_\eps^{in} (x,v) \geq 0 \qquad \ \, \qquad {\text{on}}\  \R^3_+ \times \R^3 \,,
\end{array}
\right.
\end{equation}
where the dimensionless number $\eps>0$ is the Knudsen number which is the ratio of the mean free path to the macroscopic length scale, and $F_\eps (t,x,v) \geq 0$ is the number density of particles at time $t$ with velocity $v \in \R^3$ and position $x \in \R^3_+ = \{ x \in \R^3; x_3 > 0 \}$, a half-space. The collision operator $B(F_1, F_2)$ with hard sphere collision interaction is defined by
\begin{equation}\label{Collision_Operator}
  \begin{aligned}
   B (F_1, F_2) (v) = & \iint_{\R^3 \times \mathbb{S}^2} |(v-u)\cdot \omega| F_1 (u^\prime) F_2 (v^\prime)  \d \omega \d u \\
   - & \iint_{\R^3 \times \mathbb{S}^2} |(v-u)\cdot \omega| F_1 (u) F_2 (v)  \d \omega \d u\,,
  \end{aligned}
\end{equation}
where
\begin{equation*}
  \begin{aligned}
    u^\prime = u + [(v-u) \cdot \omega] \omega\,,\ v^\prime = v + [(v-u) \cdot \omega] \omega \,, \ \omega \in \mathbb{S}^2 \,.
  \end{aligned}
\end{equation*}
The outward normal vector of $\R^3_+$ is denoted by $n= (0, 0, -1)$. The symbol $\Sigma := \p \R^3_+ \times \R^3$ is the phase space boundary of $\R^3_+ \times \R^3$, which can be split as outgoing boundary $\Sigma_+$, incoming boundary $\Sigma_-$, and grazing boundary $\Sigma_0$:
\begin{equation*}
  \begin{aligned}
   & \Sigma_+ = \{ (x, v) : x \in \p \R^3_+ \,, v \cdot n = - v_3 > 0 \} \,, \\
   & \Sigma_- = \{ (x, v) : x \in \p \R^3_+ \,, v \cdot n = - v_3 < 0 \} \,, \\
   & \Sigma_0 = \{ (x, v) : x \in \p \R^3_+ \,, v \cdot n = - v_3 = 0 \} \,.
  \end{aligned}
\end{equation*}

Let $\gamma_\pm F : = \mathbbm{1}_{\Sigma_\pm} F$. The following {\em complete diffusive boundary condition} is further imposed on the scaled equation \eqref{BE}:
\begin{equation}\label{Diffusive_BC}
\gamma_- F_\eps =  M_w (t, x, v) \int_{v \cdot n >0} \gamma_+ F_\eps ( v \cdot n ) \d v \quad {\text{on}}\ \R_+ \times \Sigma_-\,,
\end{equation}
where $M_w (t, x, v)$ is the local Maxwellian distribution function corresponding to the wall (boundary) with the form
\begin{equation}\label{M_w}
\begin{aligned}
M_w (t, x, v) = \frac{\rho_w (t, x)}{[2 \pi T_w (t, x)]^{3/2}} \exp\Big\{ - \frac{|v- u_w (t, x)|^2}{2 T_w (t, x)} \Big\} \,, \ t \in \R_+ \,, \ (x, v) \in \Sigma \,.
\end{aligned}
\end{equation}
Here $u_w = (u_{w,1}, u_{w,2}, u_{w,3}) \in \R^3$ and $ T_w$ are velocity and temperature of the boundary, respectively. The constant $\rho_w$ is such that
\begin{equation}\label{Mw-Condition}
  \begin{aligned}
    \int_{v \cdot n < 0} |v \cdot n| M_w (t,x,v) \d v \equiv 1 \,,
  \end{aligned}
\end{equation}
which means that
\begin{equation}\label{rho_w-general}
  \begin{aligned}
    \rho_w = \sqrt{\tfrac{2 \pi}{T_w}} \Big[ \exp \Big( - \tfrac{u_{w,3}}{2 T_w} \Big) + \tfrac{u_{w, 3}}{T_w} \int_{-u_{w,3}}^{+\infty} \exp\Big( -\tfrac{z^2}{2 T_w} \Big) \d z \Big]^{-1} \,.
  \end{aligned}
\end{equation}
Once
\begin{equation}\label{uw3=0}
  \begin{aligned}
    u_{w,3} = 0
  \end{aligned}
\end{equation}
is assumed, which denotes the fixed boundary wall, the $\rho_w$ in \eqref{rho_w-general} reads $\rho_w = \sqrt{\tfrac{2 \pi}{T_w}}$. In other words, under the assumption \eqref{uw3=0} for fixed boundary wall,
\begin{equation}\label{M_w2}
  \begin{aligned}
    M_w = \sqrt{\frac{2 \pi}{T_w }} \frac{1}{(2 \pi T_w)^{3/2}} \exp\Big( - \frac{|v- u_w |^2}{2 T_w } \Big) \ \textrm{ on } \R_+ \times \Sigma \,,
  \end{aligned}
\end{equation}
which is consistent with the form in Sone's book \cite{Sone-2007-Book}. It is worth noting that the condition \eqref{Mw-Condition} means that the particles absorbed by the wall will be completely released by the way of Gaussian distribution. We remark that the condition \eqref{Diffusive_BC} is a special case of the so-called Maxwell reflection condition:
\begin{equation}\label{MBC}
\gamma_- F_\eps = (1-\alpha_\eps) L \gamma_+ F_\eps + \alpha_\eps K \gamma_+ F_\eps \quad \textrm{on } \R_+ \times \Sigma_- \,,
\end{equation}
where the accommodation coefficient $\alpha_\eps \in [0,1]$ describes how much the molecules accommodate to the state of the wall, and the specular-reflection part $L \gamma_+ F_\eps$ and the diffuse-reflection part $K \gamma_+ F_\eps$ in \eqref{MBC} are
\begin{align*}
L \gamma_+ F_\eps (t, x, v) = & F_\eps (t, x, R_x v) \,, \ R_x v = v - 2 (v \cdot n) n = (\bar{v}, - v_3) \,, \\
K \gamma_+ F_\eps (t, x, v) = &  M_w (t, x, v) \int_{v \cdot n >0} \gamma_+ F_\eps ( v \cdot n ) \d v \,,
\end{align*}
respectively. The special case $\alpha_\eps=0$ corresponds to the specular reflection, while $\alpha_\eps=1$ refers to the complete diffusion, which is the case considered in this paper.

This paper addresses the limit of the scaled Boltzmann system \eqref{BE}-\eqref{Diffusive_BC} as the Knudsen number $\eps$ goes to zero. It is well-known that (see standard references for the Boltzmann equation \cite{Cercignani-1988, CIP-1994, Sone-2007-Book}, for instance), as $\eps \to 0$, the solutions $F_\eps$ of the Boltzmann equation \eqref{BE} formally converge to a local Maxwellian
\begin{equation}\label{Max-M}
  \begin{aligned}
    \M (t, x, v) = \M_{[\rho,\u, T]} (t,x,v) : = \frac{\rho(t, x)}{\big(2 \pi T(t, x)\big)^{\frac{3}{2}}}
    \exp\Big( - \frac{|v - \u (t, x)|^2}{2 T(t, x)} \Big) \,,
  \end{aligned}
\end{equation}
whose parameters $(\rho, \u, T)$ satisfy the compressible Euler system. In particular, for the domain with boundary, the key (also the difficulty) is to determine the boundary condition for the compressible Euler system. It was formally derived in Sone's book \cite{Sone-2007-Book} that for the Boltzmann equation with the complete diffusive boundary condition \eqref{Diffusive_BC}, the limiting Euler system would be imposed on the impermeable boundary condition (see also \cite{Schochet-1986-CMP}, for instance). The goal of the this paper is to rigorously justify this limit by using the method of Hilbert expansion, in the context of short time smooth solutions. We first state the main result of this paper in an informal way. The precise statement will be given in the later part of {\em Introduction}, after we introduce some required notations. 

\begin{theorem}[Informal]
	Given $(\rho^{in}, \u^{in}, T^{in}) (x)$, $(\rho_k^{in}, u_k^{in}, \theta_k^{in}) (x)$ and $(\bar{u}_k^{b,in}, \theta_k^{b, in}) (\bar{x}, \tfrac{x_3}{\sqrt{\eps}})$ $(1 \leq k \leq 4)$ sufficiently smooth, there exists a small $\eps_0 > 0$ such that when $0 < \eps < \eps_0$, for the initial data of \eqref{BE} $F_\eps^{in} (x,v) \geq 0$ with the form
	\begin{equation*}
		\begin{aligned}
			F_\eps^{in} (x,v) = \M_{[\rho^{in}, \u^{in}, T^{in}]} (x,v) + \sum_{k=1}^6 \sqrt{\eps}^k \{ F_k^{in} (x,v) + F^{b, in}_k (\bar{x}, \tfrac{x_3}{\sqrt{\eps}}, v) + F_k^{bb, in} (\bar{x}, \tfrac{x_3}{\eps}, v) \} \\
			+ \sqrt{\eps}^5 F_{R, \eps}^{in} (x,v) \,,
		\end{aligned}
	\end{equation*}
	the initial-boundary problem \eqref{BE}-\eqref{Diffusive_BC} admits a unique solution $F_\eps (t,x,v)$ with the form
	\begin{equation*}
		\begin{aligned}
			F_\eps (t, x, v) = \M_{[\rho, \u, T]} (t, x, v) + \sum_{k=1}^6 \sqrt{\eps}^k \big\{ F_k (t, x, v) + F^b_k (t, \bar{x}, \tfrac{x_3}{\sqrt{\eps}}, v) + F^{bb}_k (t, \bar{x}, \tfrac{x_3}{\eps}, v) \big\} \\
			+ \sqrt{\eps}^5 F_{R, \eps} (t, x, v) \,,
		\end{aligned}
	\end{equation*}
	where $(\rho, \u, T)$ is the smooth solution to the compressible Euler system with impermeable boundary condition and initial data $(\rho^{in}, \u^{in}, T^{in})$, i.e., \eqref{Compressible_Euler_Sys}-\eqref{BC-CEuler}-\eqref{IC-CEuler} below.
\end{theorem}

There are some {\bf remarks} on the previous informal theorem:
\begin{enumerate}
	\item The interior $F_k (t,x,v)$ and Knudsen layer $F_k^{bb} (t, \bar{x}, \tfrac{x_3}{\eps}, v)$ will be inductively determined, provided that the initial data 
	\begin{equation*}
		\begin{aligned}
			(\rho^{in}, \u^{in}, T^{in}) (x), (\rho_k^{in}, u_k^{in}, \theta_k^{in}) (x) \,, (\bar{u}_k^{b,in}, \theta_k^{b, in}) (\bar{x}, \tfrac{x_3}{\sqrt{\eps}}) 
		\end{aligned}
	\end{equation*}
	of the fluid variables of $\M_{[\rho, \u, T]}, F_k, F^b_k (1 \leq k \leq 4)$ are all given. The precise expressions will be derived later.

   \item The existence time interval of $F_\eps(t,x,v)$, say $[0, \tau]$ will be the same as that of the compressible Euler system \eqref{Compressible_Euler_Sys}-\eqref{BC-CEuler}-\eqref{IC-CEuler}. Furthermore, all the equations which solve the interior $F_k (t,x,v)$, viscous layer $F^b_k (t, \bar{x}, \frac{x_3}{\sqrt{\eps}}, v)$, Knudsen layer $F_k^{bb} (t, \bar{x}, \tfrac{x_3}{\eps}, v)$ and remainder $F_{R, \eps} (t, x, v)$ depends on the solution to this compressible Euler system. In this sense, the compressible Euler system determine the solution of the scaled Boltzmann equation.
	
	\item If the initial data of the all fluid variables are sufficiently smooth and $F_{R, \eps}^{in} (x,v)$ is uniformly bounded in small $\eps \in (0, \eps_0)$ for some suitable norms, then $\M_{[\rho, \u, T]}$, $F_k$, $F^b_k$, $F_k^{bb} (1 \leq k \leq 6)$ are all uniformly bounded in small $\eps \in (0, \eps_0)$ in some suitable functional spaces.
	
	\item The compressible Euler limit from the Boltzmann equation \eqref{BE} is understood in the following sense:
	\begin{equation*}
		\begin{aligned}
			\| F_\eps - \M_{[\rho, \u, T]} \|_X \to 0
		\end{aligned}
	\end{equation*}
    as $\eps \to 0$, where $X$ denotes a suitable functional space to be precisely stated later.
\end{enumerate}

\subsection{Compressible Euler limit from Boltzmann equation with boundary}
One of the key features of the Boltzmann equation (or more generally, kinetic equations) is their relation to the fluid equations. Generally speaking, there is an important dimensionless number for the Boltzmann equation---the Knudsen number $\eps>0$, which is the ratio of the mean free path to the macroscopic length scale, as mentioned above. The Knudsen number measures how fluid-like of the Boltzmann equation is. The smaller the Knudsen number is, the more fluid-like the Boltzmann equation is. In this sense, we call the limiting process of $\eps\rightarrow 0$ as hydrodynamic limit, or fluid limit. For the different physical scalings, the Boltzmann equation can formally converge to different fluid equations, both compressible and incompressible, or both viscous and inviscid. In other words, the solutions to the Boltzmann equation can converge to compressible/incompressible Navier-Stokes and Euler equations. For the formal derivations of these fluid equations, see \cite{BGL1}. Rigorously justifying these limiting process in different contexts of solutions (renormalized, or classical, or mild solutions) is an active research field in the last four decades.

There have been much significant progress in the limits to incompressible fluids, such as Navier-Stokes, Stokes, or even Euler. We only list some representative results. One group of results is in the framework of DiPerna-Lions renormalized solutions to the Boltzmann equation \cite{DiPerna-Lions}, i.e., the so-called BGL program, which aimed at justifying the limit to Leray solutions of the incompressible Navier-Stokes equations. This was initialized by Bardos-Golse-Levermore \cite{BGL1, BGL2}, and finished by Golse and Saint-Raymond \cite{Golse-SRM-04, Golse-SRM-09}. The corresponding results in bounded domain were carried out in \cite{Masmoudi-SRM-CPAM, Jiang-Masmoudi-CPAM}. We also should mention the incompressible Euler limits of Saint-Raymond \cite{SRM-ARMA03, SRM-IHP09}. Another group of results is in the framework of classical solutions, which are based on nonlinear energy method, semi-group method or hypercoercivity (the latter two further rely on the spectral analysis of the linearized Boltzmann operator), see \cite{Bardos-Ukai1991, Briant-JDE2015, BMM-AA2019, Gallagher-Tristani, Guo-CPAM06, Jang-Kim,JLM-CPDE, Jiang-Xu-Zhao,Levermore-Masmoudi-ARMA}.

Comparing to the incompressible limits listed above, the results to the compressible Euler system from the Boltzmann equation are much limited. This is mainly due to our still very poor understanding of the well-posedness of compressible Euler system for which we do not even know how to define weak solutions (at least for multiple spatial dimensions). The only available global-in-time solution is the celebrated BV solution of Glimm in 1965 \cite{Glimm-1965}. However, his solution was on 1-D. Regarding the higher dimensions, there are local-in-time classical solutions (see standard textbook \cite{Majda-1984}). The compressible Euler limit from the Boltzmann equation dates back to Japanese school for analytical data and Caflisch's Hilbert expansion approach which are based on the local well-posedness of compressible Euler system, see \cite{Caflish-1980-CPAM, Nishida1978}. Later improvement of Caflisch's result employing the recent progress on the $L^2\mbox{-}L^\infty$ estimate was also used to help justifying the linearized acoustic limit \cite{Guo-Jang-Jiang-2009-KRM, Guo-Jang-Jiang-2010-CPAM}. We should also mention the compressible Euler limit in the context of 1-D Riemann problem is away from the initial time in \cite{HWWY-SIMA2013}.

In the domain with boundary, the situation is much more complicated. As formally (and numerically) analyzed by Japanese school (summarized in Sone's books \cite{Sone-2002book, Sone-2007-Book}), to derive the equations of compressible Euler system from the Boltzmann equation with Maxwell reflection boundary condition by using Hilbert type expansion, in most cases, requires two coupled boundary layers: viscous and kinetic layers. The former is also called {\em Prandtl layer}, the latter is {\em Knudsen layer} in the physics literatures. The reason that these two types of boundary layers are needed can be explained as follows: as well-known, the leading order in {\em interior} is the compressible Euler system whose natural boundary condition is the impermeable condition $\u\!\cdot\!n=0$. However, the local Maxwellian governed by  Euler system with this condition does not satisfy the Maxwell reflection boundary condition, except the specular reflection case (i.e., $\alpha_\eps=0$). So, kinetic layers with thickness $\eps$ are needed. In these layers, each term satisfies the linear kinetic boundary layer equations which require {\em four} solvability conditions (by the theorem of Golse-Perthame-Sulem \cite{Golse-Perthame-Sulem-1988-ARMA}), but the number of boundary conditions for the compressible Euler system is {\em one}. This mismatch indicates there should be another layer with thickness $\sqrt{\eps}$: Prandtl layer. This name comes from the fact that in this layer, each term satisfies the famous Prandtl equations of compressible type. Specifically, if $\alpha_\eps=O(1)$ the leading term satisfies the {\em nonlinear} compressible Prandtl equations, the higher order terms satisfied the {\em linearized} compressible Prandtl equations. If $\alpha_\eps=O(\eps^\beta)$ with $\beta>0$, the boundary layers are weak, thus all the boundary terms appear in higher orders. Thus, nonlinear Prandtl equation does not appear. Recently, there are two results in this direction for rigorous proof in the context of smooth solutions, of course, in short time interval: for specular reflection $\alpha_\eps=0$, by Guo-Huang-Wang in \cite{GHW-2020}, and for small but non-zero accommodation coeffieicnt $\alpha_\eps= O(\sqrt{\eps})$, by the authors of the current paper, in \cite{JLT-arXiv-2021}.

In the current paper, we treat the case of complete diffusion boundary condition: $\alpha_\eps = 1$. The special feature of this paper is: this case is different with specular reflection \cite{GHW-2020} with $\alpha_\eps = 0$ and the case of $0< \alpha_\eps < 1$ (for example \cite{JLT-arXiv-2021}). For the case of $\alpha_\eps = 0$, the Prandtl layers are of Neumann boundary condition. For the case of $0< \alpha_\eps < 1$, the Prandtl layers are of Robin boundary condition. For the case $\alpha_\eps=1$, the Prandtl layers are of Dirichet boundary condition. More detailed explanations will be given in subsections \ref{Subsec:1.4} and \ref{Subsec:1.5} later.

\subsection{Notations}

Throughout this paper, the notation $\bar{U} = (U_1, U_2)$ is used for any vector $U = (U_1, U_2, U_3) \in \R^3$. Moreover, for simplicity of presentations, the notation $V^0 : = V |_{x_3 = 0}$ is employed for any symbol $V = V(x)$, which may be a function, vector or operator.

For the functional spaces,  $H^s$ denotes the Sobolev space $W^{s,2} (\R^3_+)$ with norm $\| \cdot \|_{H^s}$.  $\| \cdot  \|_2$ and $\| \cdot \|_\infty$ are the $L^2$-norm and $L^\infty$-norm in both $(x,v) \in \R^3_+ \times \R^3$ variables, respectively. $\l \cdot, \cdot \r$ is the $L^2$-inner product over $(x,v) \in \R^3_+ \times \R^3$. In particular, the symbol $\l \cdot, \cdot \r_v$ with subscript $v \in \R^3$ means the $L^2$-inner product in $v \in \R^3$.

In order to quantitatively describe the linear hyperbolic system \eqref{Linear-Hyperbolic-Syst} below in the half-space $\R^3_+$, let
\begin{equation*}
\p^\alpha_{t, \bar{x}} = \p_t^{\alpha_0} \p_{x_1}^{\alpha_1} \p_{x_2}^{\alpha_2} \,,
\end{equation*}
where $\alpha = (\alpha_0, \alpha_1, \alpha_2) \in \mathbb{N}^3$, and let
\begin{equation}\label{Notation_H_k}
\begin{aligned}
& \|f(t)\|^2_{\mathcal{H}^k (\R^3_+)} = \sum_{|\alpha|+i \le k} \|\p_{t, \bar{x}}^\alpha \p^i_{x_3} f(t)\|^2_{L^2(\R^3_+)}\,, \\
& \|g(t)\|^2_{\mathcal{H}^k(\R^2)} = \sum_{|\alpha|\le k} \|\p_{t, \bar{x}}^\alpha g(t)\|^2_{L^2(\R^2)}
\end{aligned}
\end{equation}
for functions $f (t) = f(t, \bar{x}, x_3)$ and $g(t) = g (t, \bar{x})$. Here $|\alpha| = \alpha_0 + \alpha_1 + \alpha_2$. Remark that, for a function $f = f(\bar{x}, x_3)$ independent of the variable $t$, $\| f \|_{\mathcal{H}^k (\R^3_+)}$ is equivalent to the standard Sobolev norm $\| f \|_{H^k (\R^3_+)}$.

While characterizing quantitatively the linear Prandtl-type system \eqref{Linear-Prandtl} below associated with the macroscopic parts of the viscous boundary layer, some new norms are required to be introduced. For $l \geq 0$,  the weighted norm is defined as
\begin{equation}\label{L2l}
	\begin{aligned}
		\| f \|^2_{L^2_l} = \int_{\R^2} \int_{\R_+} (1 + \zeta)^l | f (\bar{x}, \zeta) |^2 \d \bar{x} \d \zeta \,.
	\end{aligned}
\end{equation}
We further introduce a weighted Sobolev space $\mathbb{H}_l^r (\R^3_+)$ for any $r, l \geq 0$. Denote by the 2D multi-index $\beta = (\beta_1, \beta_2) \in \mathbb{N}^2$. For any $l, r \geq 0$, let
\begin{equation}\label{Def-l_j}
	\begin{aligned}
		l_j = l + 2 (r - j) \,, \ 0 \leq j \leq r \,.
	\end{aligned}
\end{equation}
We then introduce the norms
\begin{equation}\label{Hrl-t-3D}
	\begin{aligned}
		& \| f (t) \|^2_{l,r,n} = \sum_{2\gamma + |\beta| =r -n} \| \partial_t^\gamma \partial_{\bar{x}}^\beta \p_{\zeta}^n f (t) \|^2_{L^2_{l_r}} \ (0 \leq n \leq r) \,, \\
		& \| f (t) \|^2_{l,r} = \sum_{n=0}^r \| f (t) \|^2_{l,r,n} = \sum_{2\gamma + |\beta| + n =r} \| \partial_t^\gamma \partial_{\bar{x}}^\beta \partial_\zeta^n f (t) \|^2_{L^2_{l_r}} \,, \\
		& \| f (t) \|^2_{\mathbb{H}^r_{l, n} (\R^3_+)} = \sum_{j=0}^r \| f (t) \|^2_{l, j, n} \ (n = 0, 1, \cdots, r) \,, \\
		& \| f (t) \|^2_{\mathbb{H}^r_l (\R^3_+)} = \sum_{j=0}^r \| f (t) \|^2_{l, j} = \sum_{n=0}^r \| f (t) \|^2_{\mathbb{H}^r_{l, n} (\R^3_+)}
	\end{aligned}
\end{equation}
for function $f = f (t, \bar{x}, \zeta)$. For $g = g(\bar{x}, \zeta)$, let
\begin{equation}\label{Hrl-3D}
	\begin{aligned}
		\| g \|^2_{\mathbb{H}^r_l (\R^3_+)} = \sum_{j=0}^r \sum_{|\beta| + n = j} \| \p_{\bar{x}}^\beta \p_{\zeta}^n g \|^2_{L^2_{l_j}} \,.
	\end{aligned}
\end{equation}
Similarly, for function $h = h (t, \bar{x})$, let
\begin{equation}\label{Hr-2D}
	\begin{aligned}
		\| h (t) \|^2_{\mathbb{H}^r (\R^2)} = \sum_{j=0}^r \| h (t) \|^2_{\Gamma, j} = \sum_{j=0}^r \sum_{2 \gamma + |\beta| = j} \| \p_t^\gamma \p_{\bar{x}}^\beta h (t) \|^2_{L^2 (\R^2)} \,.
	\end{aligned}
\end{equation}

For any function $G = G(t, \bar{x}, y, v)$, with $(t, \bar{x}, y, v) \in \R_+ \times \R^2 \times \overline{\R}_+ \times \R^3$,  the Taylor expansion at $y = 0$ is
\begin{equation}
G = G^0 + \sum_{1 \leq l \leq N} \tfrac{y^l}{l !} G^{(l)} + \tfrac{y^{N+1}}{(N+1) !} \widetilde{G}^{(N+1)} \,,
\end{equation}
where the symbols
\begin{equation*}
\begin{aligned}
G^{(l)} = (\p_y^l G ) (t, \bar{x}, 0, v) \,, \ \widetilde{G}^{(N+1)} = ( \p_y^{N+1} G ) (t, \bar{x}, \eta, v) \ \textrm{for some } \eta \in (0, y).
\end{aligned}
\end{equation*}

Via the local Maxwellian $\M (t,x,v)$, the local linearized collision operator $\L$ is defined by
\begin{equation*}
\begin{aligned}
\L g = - \tfrac{1}{\sqrt{\M}} \Big\{ B (\M, \sqrt{\M} g) + B (\sqrt{\M} g, \M) \Big\} \,,
\end{aligned}
\end{equation*}
whose null space $\mathcal{N}$ is spanned by (see \cite{Caflish-1980-CPAM}, for instance)
\begin{equation*}
\begin{aligned}
\tfrac{1}{\sqrt{\rho}} \sqrt{\M} \,, \ \tfrac{v_i  - \u_i}{\sqrt{\rho T}} \sqrt{\M} \ (i=1,2,3) \,, \ \tfrac{1}{\sqrt{6 \rho}} \left\{ \tfrac{|v - \u|^2}{T} - 3 \right\} \sqrt{\M} \,.
\end{aligned}
\end{equation*}
The collision frequency $\nu (v) : = \nu (\M) (v)$ is given by
\begin{equation*}
\begin{aligned}
\nu (\M) (v) = \iint_{\R^3 \times \mathbb{S}^2} | (v - v') \cdot \omega| \M (v') \d v' \d \omega \,.
\end{aligned}
\end{equation*}
Note that
\begin{equation}\label{nu-phi}
\begin{aligned}
\nu (v) \thicksim \rho \l v \r \,,
\end{aligned}
\end{equation}
where $\l v \r = \sqrt{1 + |v|^2}$. The following weighted $L^2$-norm is introduced
\begin{equation*}
\begin{aligned}
\| g \|^2_\nu = \iint_{\R^3_+ \times \R^3} |g (x,v)|^2 \nu (v) \d x \d v \,.
\end{aligned}
\end{equation*}
Let $\mathcal{P} g$ be the $L^2_v$ projection with respect to $\mathcal{N}$. Then it is well-known (see for example \cite{Caflish-1980-CPAM}) that there exists a positive number $c_0 > 0$ such that
\begin{equation}\label{Hypocoercivity}
\begin{aligned}
\l \L g, g \r \geq c_0 \| (\mathcal{I} - \mathcal{P}) g \|^2_\nu \,,
\end{aligned}
\end{equation}
provided that $\rho$ and $T$ have positive lower and upper bounds.

\subsection{Main results}\label{Subsec:1.4}

As inspired in Caflish's work \cite{Caflish-1980-CPAM}, the Hilbert expansion approach is employed in current work. Based on the formal analyses in Section \ref{Sec:FormalAnaly} below, the complete diffusive boundary condition \eqref{Diffusive_BC} with the assumption \eqref{uw-Tw-match} will not generate the strong viscous boundary layers (nonlinear Prandtl equations). Moreover, as shown in Sone's book \cite{Sone-2007-Book}, the leading Knudsen layer vanishes, which, nevertheless, formally subjects to a nonlinear Knudsen layer equation. In other words, the so-called nonlinear Knudsen layer equation does NOT exist. Furthermore, the thickness of Knudsen boundary layers is $\eps$, but that of viscous layer is $\sqrt{\eps}$. Consequently, the solutions $F_\eps (t,x,v)$ to the scaled Boltzmann equation \eqref{BE} are therefore expanded by order $\sqrt{\eps}$. More precisely, the ansatz of the asymptotically expanded form is
\begin{equation}\label{Hilbert_Expansion1}
  \begin{aligned}
    F_\eps (t, x, v) = \M (t, x, v) + \sum_{k=1}^6 \sqrt{\eps}^k \big\{ F_k (t, x, v) + F^b_k (t, \bar{x}, \tfrac{x_3}{\sqrt{\eps}}) + F^{bb}_k (t, \bar{x}, \tfrac{x_3}{\eps}, v) \big\} \\
    + \sqrt{\eps}^5 F_{R, \eps} (t, x, v) \ge 0 \,,
  \end{aligned}
\end{equation}
where the local Maxwellian $\M (t,x,v)$ in \eqref{Max-M} for interior is governed by the compressible Euler equations
\begin{equation}\label{Compressible_Euler_Sys}
  \left\{
    \begin{aligned}
      &\p_t \rho + \div_x (\rho \u) = 0\,,\\
      &\p_t (\rho \u) + \div_x (\rho \u \otimes \u) + \nabla p = 0 \,, \\
      &\p_t \big[\rho \big( \tfrac{3}{2} T + \tfrac{1}{2} |\u|^2 \big) \big] + \div_x \big[ \rho \u \big( \tfrac{3}{2} T + \tfrac{1}{2} |\u|^2 \big) \big] + \div_x (p \u) = 0 \,, \\
      & p = \rho T
    \end{aligned}
  \right.
\end{equation}
over $(t, x) \in \R_+ \times \R^3_+$ with the impermeable boundary condition
\begin{equation}\label{BC-CEuler}
  \begin{aligned}
    \u \cdot n |_{x_3 =0} =  -\u_3 |_{x_3 = 0} = - \u_3^0 = 0 \,.
  \end{aligned}
\end{equation}
One further imposes the initial data (this can be realized by setting special form of the initial data of the Boltzmann equation \eqref{BE}):
\begin{equation}\label{IC-CEuler}
\begin{aligned}
(\rho, \u, T) (0, x) = (\rho^{in}, \u^{in}, T^{in}) (x)
\end{aligned}
\end{equation}
with compatibility condition $\u^{in} \cdot n |_{x_3} = 0$. One then can define an initial local Maxwellian
\begin{equation*}
\begin{aligned}
  \M^{in} (x,v) : = \frac{\rho^{in}(x)}{\big(2 \pi T^{in}(x)\big)^{\frac{3}{2}}}
  \exp\Big( - \frac{|v - \u^{in} (x)|^2}{2 T^{in}(x)} \Big) \,.
\end{aligned}
\end{equation*}
Denote by
\begin{equation*}
  \begin{aligned}
    \M^{in, 0} : = \M^{in} |_{x_3 = 0} \,.
  \end{aligned}
\end{equation*}
Remark that the boundary value \eqref{BC-CEuler} can eventually be derived from the solvability of nonlinear Knudsen boundary layer equations in half-space.

From \cite{Schochet-1986-CMP} or \cite{Chen-FMC-2007}, the following proposition holds.

\begin{proposition}\label{Proposition_Compressible_Euler}
	Let $s_0 \ge 3$. Assume $(\rho^{in} - \rho_\#, \u^{in}, T^{in} - T_\#) \in H^{s_0} (\R^3_+)$ satisfies $
	0 < \frac{3}{4} \rho_\# \le \rho^{in} (x) \le \frac{5}{4} \rho_\# \,,\quad
	0 < \frac{3}{4} T_\# \le T^{in} (x) \le \frac{5}{4} T_\#
	$ for some constants $\rho_\#$ and $T_\#$. Then there is a $\tau >0$ such that the compressible Euler system \eqref{Compressible_Euler_Sys}-\eqref{BC-CEuler}-\eqref{IC-CEuler} admits a unique solution
	$$
	(\rho - \rho_\#, \u, T - T_\#) \in C ( [0,\tau]; H^{s_0} (\R^3_+) ) \cap C^1 ( [0,\tau]; H^{s_0-1} (\R^3_+) )
	$$
	such that $ 0 < \tfrac{1}{2} \rho_\# \le \rho(t, x) \le \frac{3}{2} \rho_\#\,, \quad
	0 < \tfrac{1}{2} T_\# \le T(t, x) \le \frac{3}{2} T_\#
	$ hold for any $(t, x) \in [0, \tau] \times \R^3_+$. Moreover, the following estimate holds:
	\begin{align}\label{Compressible_Euler_Bound}
	E_{s_0} : = \| (\rho - \rho_\#, \u, T - T_\#) \|_{C ( [0,\tau]; H^{s_0} (\R^3_+) ) \cap C^1 ( [0,\tau]; H^{s_0-1} (\R^3_+) )} \le C_0 \,.
	\end{align}
	Here the constants $\tau$, $C_0 >0$ depend only on the $H^{s_0}$-quantity $\mathcal{E}^{in}_0 : = \| (\rho^{in} - \rho_\#, \u^{in}, T^{in} - T_\#) \|_{H^{s_0}}$.
\end{proposition}

We now outline the construction of the truncated expansion \eqref{Hilbert_Expansion1}, whose details will be given in Section \ref{Sec:FormalAnaly} later. Denote by
\begin{equation}
  \begin{aligned}
    f_i : = \frac{F_i}{\sqrt{\M}} \,, \ f^b_i : = \frac{F^b_i}{\sqrt{\M^0}} \,, \ f^{bb}_i : = \frac{F^{bb}_i}{\sqrt{\M^0}} \ (i \geq 1) \,.
  \end{aligned}
\end{equation}
In \eqref{Hilbert_Expansion1}, the higher order interior expanded terms $F_k$ $(1 \leq k \leq 6)$ satisfy
\begin{equation}\label{Relt1}
  \begin{aligned}
    f_k & = \P f_k + (\I - \P) f_k = \big[ \tfrac{\rho_k}{\rho} + u_k \cdot \tfrac{v - \u}{T} + \tfrac{\theta_k}{6 T} ( \tfrac{|v - \u|^2}{T} - 3 ) \big] \sqrt{\M} + (\I - \P) f_k \ (k = 1,2,3,4) \,, \\
    f_k & = (\I - \P) f_k \ (k =5,6) \,,
  \end{aligned}
\end{equation}
where the kinetic parts read
\begin{equation}\label{KinPart-12}
  \begin{aligned}
    (\I - \P) f_1 = & 0 \,, \\
    (\I - \P) f_2 = & \L^{-1} \Big( - \tfrac{ (\p_t + v \cdot \nabla_x) \M - B (F_1, F_1) }{ \sqrt{\M} } \Big) \,, \\
    (\I - \P) f_k = & \L^{-1} \Big( - \tfrac{(\p_t + v \cdot \nabla_x) F_{k-2} - \sum_{\substack{i+j = k\,,\\ i, j \geq 1}} B (F_i, F_j) }{\sqrt{\M}} \Big) \,, k =3,4,5,6 \,,
  \end{aligned}
\end{equation}
and fluid variables $(\rho_k, u_k, \theta_k) \in \R \times \R^3 \times \R$ of $f_k$ ($k=1,2,3,4$) obey the following linear hyperbolic system
\begin{equation}\label{Linear-Hyperbolic-Syst}
  \left\{
    \begin{aligned}
      & \p_t \rho_k + \div_x (\rho u_k + \rho_k \u) = 0 \,, \\
      & \rho \big( \p_t u_k + u_k \cdot \nabla_x \u + \u \cdot \nabla_x u_k \big) - \tfrac{\nabla_x (\rho T)}{\rho} \rho_k + \nabla_x \big( \tfrac{\rho \theta_k + 3 T \rho_k}{3} \big) = \mathcal{F}^\bot_u (f_k) \,, \\
      & \rho \big( \p_t \theta_k + \u \cdot \nabla_x \theta_k + \tfrac{2}{3} \big( \theta_k \div_x \u + 3 T \div_x u_k \big) + 3 u_k \cdot \nabla_x T \big) = \mathcal{F}^\bot_\theta (f_k)
    \end{aligned}
  \right.
\end{equation}
with boundary condition
\begin{equation}\label{BC-LinearHypo}
  \begin{aligned}
    u^0_k \cdot n = & \mathcal{J}_k (\M^0) \\
    : = & T^0 (\Psi_k + 5 T^0 \Theta_k) (t, \bar{x}, 0) - \tfrac{1}{\rho^0} \l (v \cdot n) \sqrt{\M^0} , (\mathcal{I} - \mathcal{P}^0) f_k (t, \bar{x}, 0, v)  \r_v \\
    & - \tfrac{1}{\rho^0} \int_0^\infty ( \partial_t \rho^b_{k-1} + \div_{\bar{x}} ( \rho^0 \bar{u}_{k-1}^b + \rho^b_{k-1} \bar{\u}^0 ) ) (t, \bar{x}, \zeta) \d \zeta \,,
  \end{aligned}
\end{equation}
where $(\rho, \u, T)$ is the smooth solution to the compressible Euler equations \eqref{Compressible_Euler_Sys}, $\Psi_k$, $\Theta_k$ are defined in \eqref{f_bb_k1_Coef} below, $\rho^b_{k-1}, \bar{u}^b_{k-1}$ are the viscous boundary layers of lower orders defined in Subsection \ref{Subsec:VBL} later, the source terms $\mathcal{F}^\bot_u (f_k ) \in \R^3$ and $\mathcal{F}^\bot_\theta (f_k) \in \R$ are defined as
\begin{equation}\label{Relt2}
  \begin{aligned}
    & \mathcal{F}^\bot_{u, i} (f_k) = - \sum_{j=1}^3 \p_{x_j} \l T \A_{ij}, (\I - \P) f_k \r_v \ (i = 1,2,3) \,, \\
    & \mathcal{F}^\bot_\theta (f_k) = - \div_x \big( 2 T^{\frac{3}{2}} \l \B, (\I - \P) f_k \r_v + 2 T \u \cdot \l \A , (\I - \P) f_k \r_v \big) - 2 \u \cdot \mathcal{F}^\bot_u (f_k)\,.
\end{aligned}
\end{equation}
Here $\A \in \R^{3 \times 3}$ and $\B \in \R^3$ are the Burnett functions with entries
\begin{equation}
  \begin{aligned}
    \A_{ij} = & \Big\{ \tfrac{(v_i - \u_i) (v_j - \u_j)}{T} - \delta_{ij} \tfrac{|v - \u|^2}{3 T} \Big\} \sqrt{\M} \quad (1 \leq i, j \leq 3) \,, \\
    \B_i = & \tfrac{v_i - \u_i}{2 \sqrt{T}} \Big( \tfrac{|v - \u|^2}{T} - 5 \Big) \sqrt{\M} \quad (1 \leq i \leq 3) \,.
  \end{aligned}
\end{equation}
They are both in $\mathcal{N}^\perp$ and perpendicular to each other. Remark that the boundary condition \eqref{BC-LinearHypo} is derived from the solvability of the higher order Knudsen layer equation $F^{bb}_k$ ($k=1,2,3,4,5$), see Section \ref{Sec:FormalAnaly} below. We also emphasize that the decay estimates for the inverse operator $\L^{-1}$ in \eqref{KinPart-12} can be found in \cite{JLT-arXiv-2021}, which is important to controlling the expanded terms. Finally, for $k = 1,2,3,4$, the initial data of \eqref{Linear-Hyperbolic-Syst} are imposed on
\begin{equation}\label{IC-Linear-Hyperbolic}
  \begin{aligned}
    (\rho_k, u_k, \theta_k) (0, x) = (\rho_k^{in}, u_k^{in}, \theta_k^{in}) (x) \in \R \times \R^3 \times \R
  \end{aligned}
\end{equation}
with compatibility $u_k^{in} \cdot n |_{x_3 = 0} = \mathcal{J}_k (\M^{in, 0})$, where the symbols $\mathcal{J}_k$ are given in \eqref{BC-LinearHypo}.

The viscous boundary layers $f^b_k (t, \bar{x}, \zeta, v)$ with the normal coordinate $\zeta = \frac{x_3}{\sqrt{\eps}}$ can be split by
\begin{equation*}
	\begin{aligned}
		f^b_k = & \left\{ \tfrac{\rho^b_k}{\rho^0} + u^b_k \cdot \tfrac{v - \u^0}{T^0} + \tfrac{\theta^b_k}{6 T^0} ( \tfrac{|v - \u^0|^2}{T^0} - 3 ) \right\} \sqrt{\M^0} + (\I - \P^0) f^b_k \,, k =1,2,3,4 \,, \\
		f^b_k = & (\I - \P^0) f^b_k \,, k = 4,5 \,.
	\end{aligned}
\end{equation*}
Denote by
\begin{equation*}
	\begin{aligned}
		p^b_k = \tfrac{\rho^0 \theta^b_k + 3 T^0 \rho^b_k}{3} \,.
	\end{aligned}
\end{equation*}
Then the kinetic parts read
\begin{equation*}
	\begin{aligned}
		(\I - \P^0) f^b_1 = & 0 \,, \\
		(\I - \P^0) f^b_k = & (\mathcal{L}^0)^{-1} \Big\{ - (\I - \P^0) (v_3 \p_\zeta \P^0 f^b_{k-1}) + \tfrac{\zeta}{\sqrt{\M^0}} \big[ B(\M^{(1)}, \sqrt{\M^0} \P^0 f^b_{k-1}) \\
		& + B(\sqrt{\M^0} \P^0 f^b_{k-1}, \M^{(1)})\big] + \tfrac{1}{\sqrt{\M^0}} \big[ B(F^0_1, \sqrt{\M^0} \P^0 f^b_{k-1}) + B(\sqrt{\M^0} \P^0 f^b_{k-1}, F^0_1) \big] \\
		& + \tfrac{1}{\sqrt{\M^0}} \big[ B(\sqrt{\M^0} f^b_1, \sqrt{\M^0} \P^0 f^b_{k-1}) + B(\sqrt{\M^0} \P^0 f^b_{k-1}, \sqrt{\M^0} f^b_1) \big] \Big\} + J^b_{k-2}\,,
	\end{aligned}
\end{equation*}
with $k = 2,3,4,5,6$, where the quantity $J^b_{k-2}$ is defined in \eqref{J_b_k-1} later. $u^b_{k,3}$ and $p^b_k$ ($k = 1,2,3,4$) are determined by 
\begin{equation*}
	\begin{aligned}
		& u^b_{1,3} = p^b_1 = 0 \,, \\
		& \p_\zeta u^b_{k, 3} = - \frac{1}{\rho^0} ( \p_t \rho^b_{k-1} + \div_{\bar{x}} ( \rho^0 \bar{u}^b_{k-1} + \rho^b_{k-1} \bar{\u}^0 ) )\,, \ \lim_{\zeta \to \infty} u^b_{k,3} (t, \bar{x}, \zeta) = 0 \,,
	\end{aligned}
\end{equation*}
and
\begin{equation*}
	\left\{
	\begin{aligned}
		\p_\zeta p^b_k
		= - \rho^0 \p_t u^b_{k-1,3} - \rho^0 \bar{\u}^0 \!\cdot\! \nabla_{\bar{x}} u^b_{k-1,3} + \rho^0 \p_{x_3} \u^0_3 u^b_{k-1,3} + \tfrac{4}{3} \mu(T^0) \p_{\zeta\zeta} u^b_{k-1,3} \\
		- \tfrac{4}{3} \rho^0 \p_\zeta \Big[ \big(\p_{x_3} \u^0 \zeta + u^0_{1,3} \big) u^b_{k-1,3} \Big] - \p_\zeta \l T^0 \A^0_{33}, J^b_{k-2} \r_v + W^b_{k-2,3} \,, \\
		\lim_{\zeta \to \infty} \theta^b_k (t, \bar{x}, \zeta) = 0 \,,
	\end{aligned}
	\right.
\end{equation*}
Here the quantities $W^b_{k-2,3}$ and $J^b_{k-2}$ are defined in \eqref{W_b_k-1} and \eqref{J_b_k-1} below. Moreover, the variables $(u^b_{k,1}, u^b_{k,2}, \theta^b_k)$ $(k = 1, 2, 3, 4)$ satisfy the following {\em linear compressible Prandtl-type} equations
\begin{equation*}
	\left\{
	\begin{aligned}
		\rho^0 (\p_t + \bar{\u}^0 \cdot \nabla_{\bar{x}}) u^b_{k,i} + \rho^0 (\p_{x_3} \u^0_3 \zeta + u^0_{1,3}) \p_{\zeta} u^b_{k,i} + & \rho^0 \bar{u}^b_k \cdot \nabla_{\bar{x}} \u^0_i + \tfrac{\p_{x_3} p^0}{3 T^0} \theta^b_k \\
		= & \mu (T^0) \p_{\zeta}^2 u^b_{k,i} + \mathsf{f}^b_{k-1, i} \ (i=1,2) \,, \\
		\rho^0 \p_t \theta^b_k + \rho^0 \bar{\u}^0 \cdot \nabla_{\bar{x}} \theta^b_k + \rho^0 \big( \p_{x_3} \u^0 \zeta + \u^0_{1,3} \big) \p_\zeta \theta^b_k + & \tfrac{2}{3} \rho^0 \div_{x} \u^0 \theta^b_k \\
		= & \tfrac{3}{5} \kappa(T^0) \p_{\zeta\zeta} \theta^b_k + \mathsf{g}^b_{k-1} \,, \\
		\lim_{\zeta \to \infty} (\bar{u}^b_k, \theta^b_k) (t, \bar{x}, \zeta) = & 0 \,,
	\end{aligned}
	\right.
\end{equation*}
with the boundary condition \eqref{BC-Prandtl} and initial data \eqref{IC-Prandtl}, i.e.,
\begin{equation*}
	\begin{aligned}
		( \bar{u}^b_k, \theta^b_k ) (t, \bar{x}, 0) = ( \bar{u}_w, T_w ) (t, \bar{x}) \,, \quad (\bar{u}^b_k, \theta^b_k) (0, \bar{x}, \zeta) = (\bar{u}^{b, in}_k , \theta^{b, in}_k ) (\bar{x}, \zeta) \,.
	\end{aligned}
\end{equation*}
Here the quantities $\mathsf{f}^b_{k-1, i}$ and $\mathsf{g}^b_{k-1}$ are given in \eqref{fg-k} below.

The Knudsen layers $f^{bb}_k (t, \bar{x}, \xi, v)$ with normal coordinate $\xi = \tfrac{x_3}{\eps}$ can be decomposed by
\begin{equation}\label{Relt3}
  \begin{aligned}
    f^{bb}_k = f^{bb}_{k,1} + f^{bb}_{k,2} \ (1 \leq k \leq 6) \,,
  \end{aligned}
\end{equation}
where
\begin{equation}\label{Relt4}
  \begin{aligned}
    f^{bb}_{k,1} = \big\{ \Psi_k v_3 + \Phi_{k,1} v_3 (v_1 - \u^0_1) + \Phi_{k,2} v_3 (v_2 - \u^0_2) + \Phi_{k,3} + \Theta_k v_3 |v- \u^0|^2 \big\} \sqrt{\M^0} \,,
  \end{aligned}
\end{equation}
and $f^{bb}_{k,2}$ $(1 \leq k \leq 6)$ subject to the Knudsen boundary layer equation
\begin{equation}\label{f_bb_k_2}
  \left\{
    \begin{aligned}
      & v_3 \p_\xi f^{bb}_{k, 2} + \mathcal{L}^0 f^{bb}_{k, 2} = S^{bb}_{k, 2} - \big( v_3 \p_\xi f^{bb}_{k,1} - S^{bb}_{k, 1} + \mathcal{L}^0 f^{bb}_{k,1} \big) \in (\mathcal{N}^0)^\perp \,, \\
      & f^{bb}_{k,2} (t, \bar{x}, 0, v) |_{v_3 >0} = \mathcal{D}_w f^{bb}_{k, 2} (t, \bar{x}, 0, v) + \mathbbm{f}_k (t, \bar{x}, v) \,, \\
      & \lim_{\xi \to +\infty} f^{bb}_{k,2} (t, \bar{x}, \xi, v) =0 \,.
    \end{aligned}
  \right.
\end{equation}
Here $\Psi_k$, $\Phi_{k,1}$, $\Phi_{k,2}$, $\Phi_{k,3}$, $\Theta_k$ are given in \eqref{f_bb_k1_Coef} below, which depend on the quantity $S^{bb}_{k,1} \in \mathcal{N}^0$ (see \eqref{S_bb_k+1_1+2} below). $S^{bb}_{k,2}$ is defined in \eqref{S_bb_k+1_1+2} below with $N=5$. The diffusive boundary operator $\mathcal{D}_w$ is defined by
\begin{equation}\label{Dw}
  \begin{aligned}
    \mathcal{D}_w f = \tfrac{M_w}{\sqrt{\M^0}} \int_{v' \cdot n >0} (v' \cdot n) f \sqrt{\M^0} \d v' \,,
  \end{aligned}
\end{equation}
where $M_w$ is the wall Maxwellian distribution. The boundary source term $\mathbbm{f}_k (t, \bar{x}, v)$ is
\begin{equation}\label{fk}
  \begin{aligned}
    \mathbbm{f}_k (t, \bar{x}, v) = - ( \I - \mathcal{D}_w ) ( f_k + f^b_k + f^{bb}_{k,1} ) (t, \bar{x}, 0, v) \,.
  \end{aligned}
\end{equation}

The remainder $F_{R, \eps} (t,x,v)$ obeys
\begin{equation}\label{Remainder_Eq}
  \begin{aligned}
    & \p_t F_{R, \eps} + v \cdot \nabla_x F_{R, \eps} - \tfrac{1}{\eps} \big[ B(F_{R, \eps}, \M) + B(\M, F_{R, \eps}) \big] \\
    = & \sum_{i=1}^6 \sqrt{\eps}^{i-2} \big[ B(F_{R, \eps}, F_i + F^{bb}_i) + B(F_i + F^{bb}_i, F_{R, \eps}) \big] \\
    & + \sqrt{\eps}^3 B(F_{R, \eps}, F_{R, \eps}) + R_\eps + R_\eps^b + + R^{bb}_\eps
  \end{aligned}
\end{equation}
with the complete diffusive boundary condition
\begin{equation}\label{Remainder_BC}
  \gamma_- F_{R, \eps} = M_w \int_{v^\prime \cdot n >0} (v^\prime \cdot n) \gamma_+ F_{R, \eps} \d v^\prime\,\  {\text{on}}\ \Sigma_-\,,
\end{equation}
where
\begin{equation}\label{R_eps}
  R_\eps = - ( \p_t + v \cdot \nabla_x ) ( F_5 + \sqrt{\eps} F_6 ) + \sum_{\substack{ i + j \geq 7, \\ 1 \leq i,j \leq 6}} B(F_i, F_j) \,,
\end{equation}
and
\begin{equation}\label{R_eps_b}
	\begin{aligned}
		R^b_\eps = & - (\p_t + \bar{v} \cdot \nabla_{\bar{x}}) (F^b_5 + \sqrt{\eps} F^b_6) - v_3 \p_\zeta F^b_6 \\
		& + \sum_{\substack{ j+l \geq 7 \\ 1 \leq j \leq 6, 1 \leq l \leq 5 }} \sqrt{\eps}^{j+l-7} \tfrac{\zeta^l}{l !} [B (\M^{(l)}, F^b_j) + B(F^b_j, \M^{(l)})] \\
		& + \sum_{\substack{ i+j \geq 7 \\ 1 \leq i, j \leq 6 }} \sqrt{\eps}^{i+j-7} [ B(F^0_i, F^b_j) + B(F^b_j, F^0_i) + B(F^b_i, F^b_j) ] \\
		& + \sum_{\substack{ i+j + l \geq 7 \\ 1 \leq i, j \leq 6, 1 \leq l \leq 5 }} \sqrt{\eps}^{i+j+l-7} \tfrac{\zeta^l}{l !} [B(F_i^{(l)}, F^b_j) + B(F^b_j, F_i^{(l)})] \\
		& + \tfrac{\zeta^6}{6!} \sum_{j=1}^6 \sqrt{\eps}^{j-1} [ B (\widetilde{\M}^{(6)} + \sum_{i=1}^6 \sqrt{\eps}^i \widetilde{F}^{(6)}_i , F^b_j) + B (F^b_j, \widetilde{\M}^{(6)} + \sum_{i=1}^6 \sqrt{\eps}^i \widetilde{F}^{(6)}_i) ] \,,
	\end{aligned}
\end{equation}
and
\begin{align}\label{R_eps_bb}
  \no R^{bb}_\eps = & - (\p_t + \bar{v} \cdot \nabla_{\bar{x}}) (F^{bb}_5 + \sqrt{\eps} F^{bb}_6) \\
  \no & + \sum_{\substack{ j+2l \geq 7 \\ 1 \leq j \leq 6, 1 \leq l \leq 5 }} \sqrt{\eps}^{j+2l-7} \tfrac{\xi^l}{l!} [ B(\M^{(l)}, F^{bb}_j) + B(F^{bb}_j, \M^{(l)}) ] \\
  \no & + \sum_{\substack{ i+j \geq 7 \\ 1 \leq i, j \leq 6 }} \sqrt{\eps}^{i+j-7} [B(F^0_i + F^{b,0}_i, F^{bb}_j) + B(F^{bb}_j, F^0_i + F^{b,0}_i) + B( F^{bb}_i, F^{bb}_j ) ] \\
  \no & + \sum_{\substack{ i+j+2l \geq 7 \\ 1 \leq i, j \leq 6, 1 \leq l \leq 5 }} \sqrt{\eps}^{i+j+2l-7} \tfrac{\xi^l}{l!} [ B(\widetilde{F}^{(l)}_i, F^{bb}_j) + B (F^{bb}_j, \widetilde{F}^{(l)}_i) ] \\
  \no & + \sum_{\substack{ i+j+l \geq 7 \\ 1 \leq i, j \leq 6, 1 \leq l \leq 5 }} \sqrt{\eps}^{i+j+l-7} \tfrac{\xi^l}{l!} [ B(\widetilde{F}^{b, (l)}_i, F^{bb}_j) + B (F^{bb}_j, \widetilde{F}^{b, (l)}_i) ] \\
  \no & + \tfrac{\xi^6}{6!} \sum_{j=1}^6 \sqrt{\eps}^{j-1} [ B  (\sqrt{\eps}^6 \widetilde{\M}^{(6)} + \sum_{i=1}^6 \sqrt{\eps}^{i+6} \widetilde{F}^{(6)}_i + \sqrt{\eps}^i \widetilde{F}^{b, (6)}_i , F^{bb}_j) \\
  & \qquad \qquad \qquad \qquad + B (F^{bb}_j, \sqrt{\eps}^6 \widetilde{\M}^{(6)} + \sum_{i=1}^6 \sqrt{\eps}^{i+6} \widetilde{F}^{(6)}_i + \sqrt{\eps}^i \widetilde{F}^{b, (6)}_i) ] \,.
\end{align}
Furthermore, the following initial data are imposed on the remainder equation \eqref{Remainder_Eq}:
\begin{equation}\label{IC-Remainder-F}
  \begin{aligned}
    F_{R, \eps} (0, x, v) = F_{R, \eps}^{in} (x,v) \,,
  \end{aligned}
\end{equation}
which satisfy the compatibility condition
\begin{equation*}
  \begin{aligned}
    \gamma_- F_{R, \eps}^{in} = M_w \int_{v^\prime \cdot n >0} (v^\prime \cdot n) \gamma_+ F_{R, \eps}^{in} \d v^\prime
  \end{aligned}
\end{equation*}
on $\Sigma_-$.

For the remainder $F_{R, \eps}$, let
\begin{equation}\label{f_Reps-h_Reps}
  \begin{aligned}
    f_{R, \eps} = \tfrac{F_{R, \eps}}{\sqrt{\M}} \,, \quad h_{R, \eps}^\ell = \l v \r^\ell \tfrac{F_{R,\eps}}{\sqrt{\M_M}}
  \end{aligned}
\end{equation}
for $\ell \geq 7$, where the global Maxwellian $\M_M = \M_M (v)$ is introduced by \cite{Caflish-1980-CPAM}
\begin{align}\label{Global_Maxwellian_M}
  \M_M = \tfrac{1}{(2\pi T_M)^{\frac{3}{2}}} \exp \Big\{ - \tfrac{|v|^2}{2 T_M} \Big\} \,.
\end{align}
Here the positive constant $T_M$ satisfies
\begin{align} \label{T_M}
  T_M < \max_{t \in [0, \tau], x \in \R^3_+} T (t, x) < 2 T_M \,.
\end{align}
Then there exist positive constants $C_1$, $C_2$ such that for some $\frac{1}{2} < z <1$ and each $(t, x, v) \in [0, \tau] \times \R^3_+ \times \R^3$, the following inequality holds:
\begin{align}\label{M-Bound}
  C_1 \M_M \le \M \le C_2 ( \M_M )^z\,.
\end{align}

One then introduces some fixed functions $F^{in}_k (x,v)$, $F^{b, in}_k (\bar{x}, \zeta, v)$ and $F^{bb, in}_k (\bar{x}, \xi, v)$ ($1 \leq k \leq 6$) constructed by the initial data $(\rho^{in}, \u^{in}, T^{in})$, $(\rho_k^{in}, u_k^{in}, \theta_k^{in})$ and $(\bar{u}^{b, in}_k , \theta^{b, in}_k ) $ ($k=1,2,3,4$) mentioned above. Specifically, the construction of $F^{in}_k (x,v)$, $F^{b, in}_k (\bar{x}, \zeta, v)$ and $F^{bb, in}_k (\bar{x}, \xi, v)$ is the same as that of $F_k (t,x,v)$, $F^b_k (t, \bar{x}, \zeta, v)$ and $F^{bb}_k (t,\bar{x}, \xi, v)$ above. It is just to replace $(\rho, \u, T)$, $(\rho_k, u_k, \theta_k)$ and $ ( \bar{u}^b_k, \theta^b_k )$ with $(\rho^{in}, \u^{in}, T^{in})$, $(\rho_k^{in}, u_k^{in}, \theta_k^{in})$ and $(\bar{u}^{b, in}_k , \theta^{b, in}_k )$, respectively. Then the initial data $F^{in}_\eps (x,v)$ of \eqref{BE} is assumed to be well-prepared:
\begin{equation}\label{IC-wellprepared}
  \begin{aligned}
    F^{in}_\eps (x,v) = \M^{in} (x, v) + \sum_{k=1}^6 \sqrt{\eps}^k \big\{ F_k^{in} (x, v) + F^{b, in}_k (\bar{x}, \tfrac{x_3}{\sqrt{\eps}}, v) + F^{bb, in}_k (\bar{x}, \tfrac{x_3}{\eps}, v) \big\} \\
    + \sqrt{\eps}^5 F_{R, \eps}^{in} (x, v) \ge 0 \,,
  \end{aligned}
\end{equation}
Furthermore, the wall velocity $u_w$ and temperature $T_w$ are assumed to match the interior fluid velocity $\u$ and temperature $T$ constrained to the boundary $\p \R^3_+$, i.e.,\footnote{ In fact, one can also assume $(u_w, T_w)$ to be a perturbation around $(\u^0, T^0)$, in which all the processes of proof are the same except that one should deal with the additional small perturbed quantities.}
\begin{equation}\label{uw-Tw-match}
  \begin{aligned}
    u_w = \u^0 \,, \ T_w = T^0  \,,
  \end{aligned}
\end{equation}
which is such that the strong boundary layers vanished. For more details, see Subsection \ref{Subsec:2.2} below.

Now the main theorem is stated as follows.

\begin{theorem}\label{Main-Thm}
	Consider the wall Maxwellian $M_w$ given in \eqref{M_w2} and \eqref{uw-Tw-match}. Assume $\ell \geq 7$. Let $(\rho^{in}, \u^{in}, T^{in})$ satisfies the assumptions in Proposition \ref{Proposition_Compressible_Euler}, and let
	\begin{equation*}
	  \begin{aligned}
	    \mathcal{E}_1^{in} : = \sum_{j=1}^4 \big( \| (\rho_j^{in}, u_j^{in}, \theta_j^{in}) \|_{\mathcal{H}^{s_j} (\R^3_+) } + \| (\bar{u}^{b, in}_j, \theta_j^{b, in}) \|_{\mathbb{H}^{s_j^b}_{l^b_j} (\R^3_+)} \big) < \infty 
	  \end{aligned}
	\end{equation*}
	for some sufficiently large numbers $s_j, s^b_j, l^b_j$ given in Proposition \ref{Prop-HigherOrd}. Correspondingly,
	\begin{itemize}
		\item $(\rho, \u, T) (t,x)$ is the solution to the compressible Euler system \eqref{Compressible_Euler_Sys}-\eqref{BC-CEuler}-\eqref{IC-CEuler} over time interval $t \in [0, \tau]$ constructed in Proposition \ref{Proposition_Compressible_Euler}, which governs the local Maxwellian $\M (t,x,v)$ in \eqref{Max-M};
		
		\item $F_k (t, x, v)$, $F^b_k (t, \bar{x}, \tfrac{x_3}{\sqrt{\eps}}, v)$ and $F^{bb}_k (t, \bar{x}, \tfrac{x_3}{\eps}, v)$ $(1 \leq k \leq 6)$ are constructed in Proposition \ref{Prop-HigherOrd} below.
	\end{itemize}
    Then there is a small $\eps_0 > 0$ such that if
    \begin{equation*}
      \begin{aligned}
        \mathcal{E}_R^{in} : = \sup_{\eps \in (0, \eps_0)} \Big( \| \tfrac{F_{R,\eps}^{in}}{\sqrt{\M^{in}}} \|_2 + \sqrt{\eps}^3 \| \l v \r^\ell \tfrac{F_{R,\eps}^{in}}{\sqrt{\M^{in}}} \|_\infty \Big) < \infty \,,
      \end{aligned}
    \end{equation*}
    the scaled Boltzmann equation \eqref{BE} with complete diffusive boundary condition \eqref{Diffusive_BC} and well-prepared initial data \eqref{IC-wellprepared} admits a unique solution $F_\eps (t,x,v)$ for $\eps \in (0, \eps_0)$ over the time interval $t \in [0, \tau]$ with the expanded form \eqref{Hilbert_Expansion1}, i.e.,
    \begin{equation*}
      \begin{aligned}
        F_\eps (t, x, v) = \M (t, x, v) + \sum_{k=1}^6 \sqrt{\eps}^k \big\{ F_k (t, x, v) + F^b_k ( t, \bar{x}, \tfrac{x_3}{\sqrt{\eps}}, v ) + F^{bb}_k (t, \bar{x}, \tfrac{x_3}{\eps}, v) \big\} \\
        + \sqrt{\eps}^5 F_{R, \eps} (t, x, v) \ge 0 \,,
      \end{aligned}
    \end{equation*}
    where the remainder $F_{R,\eps} (t,x,v)$ satisfies the uniform bound
    \begin{equation*}
      \begin{aligned}
        \sup_{t \in [0, \tau]} \Big\{ \| \tfrac{F_{R,\eps}}{\sqrt{\M}} (t) \|_2 + \sqrt{\eps}^3 \| \l v \r^\ell \tfrac{F_{R,\eps}}{\sqrt{\M_M}} (t) \|_\infty \Big\} \leq C( \tau, \mathcal{E}_0^{in}, \mathcal{E}_1^{in}, \mathcal{E}_R^{in} ) < \infty \,.
      \end{aligned}
    \end{equation*}
    Here the quantity $\mathcal{E}_0^{in}$ is given in Proposition \ref{Proposition_Compressible_Euler}.
\end{theorem}

\begin{remark}
	Together with Proposition \ref{Prop-HigherOrd} below, Theorem \ref{Main-Thm} shows that as $\eps \to 0$,
	\begin{equation*}
	  \begin{aligned}
	    \sup_{t \in [0, \tau]} \Big\{ \| (\tfrac{F_\eps - \M}{\sqrt{\M}}) (t) \|_2 + \| \l v \r^\ell ( \tfrac{F_\eps - \M}{\sqrt{\M_M}} ) (t) \|_\infty \Big\} \leq C {\eps} \to 0 \,.
	  \end{aligned}
	\end{equation*}
	Therefore, while imposed on the well-prepared initial data, one has justified the hydrodynamic limit from the Boltzmann equation with complete diffusive boundary condition to the compressible Euler system with impermeable boundary condition for the half-space problem. Moreover, the convergence rate is at least $\mathcal{O} (\eps )$.
\end{remark}

\subsection{Sketch of ideas}\label{Subsec:1.5}

As indicated in Sone's book \cite{Sone-2007-Book}, there is no leading Knudsen layer for the Boltzmann equation with the complete diffusive boundary condition in half-space. Actually, even if the leading Knudsen layer is supposed to exist, the one formally obeys a nonlinear equation. However, it can be proved that the solution must be zero\footnote{This is the reason why the nonlinear Knudsen boundary layer problem under the diffusive reflection boundary condition in half-space is not needed. We give an explanation in Appendix \ref{Appendix} that even if we ``pretend" to assume the $O(1)$ Knudsen layer exists, the corresponding nonlinear kinetic boundary layer equation would have only zero solution.}. So, the Knudsen layers are initially expanded from the next higher order, which satisfy the linear equations with the diffusive boundary condition. Furthermore, under the assumption \eqref{uw-Tw-match}, the strong viscous boundary layer do NOT exist, see Subsection \ref{Subsec:2.2}. But the weak viscous boundary layers at higher orders are required, whose thickness is $\sqrt{\eps}$. Note that the thickness of the Knudsen layer is $\eps$. Taking the two facts above into consideration, the thickness of expansions is $\sqrt{\eps}$. Moreover, there only requires one solvability condition for the linear Knudsen layer equations with diffusive boundary condition, see \cite{Coron-Golse-Sulem-1988-CPAM}. The single solvability condition will determine one boundary condition such that the interior expansion is solvable. Then, the Prandtl layers will be artificially imposed on the Dirichet boundary values such that the complete diffusive reflection boundary condition holds at every expanded order. Therefore, the ansatz of the expanded form \eqref{Hilbert_Expansion1} is done.


In $L^2$ estimates for the remainder equation, the key point is to deal with the boundary integral $- \tfrac{1}{2} \iint_\Sigma v_3 |f_{R, \eps} (t,\bar{x}, 0,v)|^2 \d v \d \bar{x}$. Thanks to the conditions \eqref{Mw-Condition} and \eqref{uw-Tw-match} satisfied by the wall Maxwellian $M_w (t, \bar{x}, v)$, namely, $\int_{v \cdot n < 0} |v \cdot n| M_w \d v = 1$ and $(u_w, T_w) = (\u^0, T^0)$, the boundary integral plays a role of boundary dissipative rate
\begin{equation*}
  \begin{aligned}
    \tfrac{1}{2} \iint_{\Sigma_+} |v \cdot n| |(\I - \mathcal{D}_w) f_{R,\eps} (t, \bar{x}, 0, v)|^2 \d v \d \bar{x} \,.
  \end{aligned}
\end{equation*}

In $L^\infty$ estimates for the remainder equation, the goal is to derive the pointwise bounds for $h_{R,\eps}^\ell (t,x,v)$ on $(t,x,v) \in [0, \tau] \times \R^3_+ \times \R^3$ along the trajectories over $s \in [0,t]$. The set $ [0, \tau] \times \R^3_+ \times \R^3$ is split into two parts: $\mathscr{U}$ and $\mathscr{V}$, see \eqref{UV} below. In fact, as time evolves, the points of the set $\mathscr{U}$ are almost all transformed into that of the set $\mathscr{V}$ through the boundary effect. For fixed $(t,x,v) \in \mathscr{U}$, the trajectory on $[0, t]$ does not hit the boundary $\p \R^3_+$. The trajectory on $[0,t]$ is a straight line segment. Following the almost same arguments in \cite{Guo-Jang-Jiang-2010-CPAM}, the required estimates on $\mathscr{U}$ can be obtained. For fixed $(t,x,v) \in \mathscr{V}$, the trajectory will hit the boundary $\p \R^3_+$ at time $t_1 = t - t_b (x,v) > 0$. The way is to integrate the $h_{R,\eps}^\ell$ equation along the trajectory on $[t_1, t]$, which is still a straight line segment. In this case, one should focus on the term $\exp\Big\{ -\tfrac{1}{\eps} \int_{t_1}^t \nu (\phi) \d \phi \Big\} h^\ell_{R, \eps} (t_1, x_b (x,v), v)$, where $(x_b (x,v), v) \in \Sigma_-$. Together with the diffusive boundary condition, $h^\ell_{R, \eps} (t_1, x_b (x,v), v)$ converts to the integral form $M_w \tfrac{\l v \r^\ell}{\sqrt{\M_M}} \int_{v^\prime \cdot n >0} (v^\prime \cdot n) \tfrac{\sqrt{\M_M}}{\l v^\prime \r^\ell} h^\ell_{R, \eps} (t_1, x_b (x,v), v^\prime) \d v^\prime$. Observe that $(t_1, x_b (x,v), v^\prime) \in \mathscr{U}$. Therefore the estimates on $\mathscr{U}$ can be applied to control this quantity. Collecting the estimates on $\mathscr{U}$ and $\mathscr{V}$, one finishes the $L^\infty$ estimates.

\subsection{Organizations of current paper}

The next section focuses on the formal analysis, in which the impermeable boundary condition of the compressible Euler system and the boundary conditions of the fluid variables of higher orders are derived. The corresponding truncation is also done. In Section \ref{Sec:F1F2}, the expanded terms $F_k$, $F^b_k$ and $F^{bb}_k$ ($1 \leq k \leq 6$) are controlled in some sufficiently smooth norms, the fact of which is such that the coefficients consisted of $F_k$, $F^b_k$ and $F^{bb}_k$ ($1 \leq k \leq 6$) appeared in the remainder equation can be bounded by some harmless constant $C > 0$. In Section \ref{Sec:Remainder}, the bounds of remainder $F_{R,\eps}$ (uniform in small Knudsen number $\eps$) are obtained in the $L^2$-$L^\infty$ framework.

\section{Formal Analysis}\label{Sec:FormalAnaly}

In this section, the goal is to derive the formal analysis in the Hilbert expansion framework. Notice that the thickness of the Knudsen layer is $\eps$, but that of the viscous layer is $\sqrt{\eps}$. As shown in Sone's book \cite{Sone-2007-Book}, there is not leading Knudsen boundary layer, which means the layers start from the order of $\sqrt{\eps}$. Moreover, the assumption \eqref{uw-Tw-match} is such that the strong viscous layer vanishes. Then the viscous layers also start from the order of $\sqrt{\eps}$. As a result, the asymptotically expanded form will be
\begin{equation}\label{Exp1}
  \begin{aligned}
    F_\eps (t, x, v) = F_0 (t,x,v) + \sum_{k \geq 1} \sqrt{\eps}^k \big[ F_k (t,x,v) + F^b (t, \bar{x}, \zeta, v) + F^{bb}_k (t, \bar{x}, \xi, v) \big]
  \end{aligned}
\end{equation}
with the normal coordinate
\begin{equation*}
  \begin{aligned}
   \zeta = \tfrac{x_3}{\sqrt{\eps}} \,, \ \xi = \tfrac{x_3}{\eps} \,,
  \end{aligned}
\end{equation*}
where $F_k$, $F^b_k$ and $F^{bb}_k$ are the interior, viscous boundary layer (also Prandtl boundary layer) and Knudsen boundary layers expansions, respectively.

\subsection{Interior expansion}

From plugging the interior expansion 
\begin{equation*}
	\begin{aligned}
		F_\eps (t,x,v) \thicksim \sum_{k \ge 0} \sqrt{\eps}^k F_k (t, x, v)
	\end{aligned}
\end{equation*}
into \eqref{BE} and collecting the same orders,
  \begin{align}\label{Interior_Ord_Ana}
     \no \sqrt{\eps}^{-2}: \qquad & 0 = B(F_0, F_0) \,,\\
     \no \sqrt{\eps}^{-1}: \qquad & 0 = B (F_0, F_1) + B (F_1, F_0) \,, \\
     \no \sqrt{\eps}^{0}: \qquad & (\p_t + v \cdot \nabla_x) F_0 = B(F_0, F_2) + B(F_2, F_0)\,,\\
     \sqrt{\eps}^{1}: \qquad & (\p_t + v \cdot \nabla_x) F_1 = B(F_0, F_3) + B(F_3, F_0) + B(F_1, F_2) + B(F_2, F_1)\,,\\
     \no ......&\\
     \no \sqrt{\eps}^k: \qquad &  (\p_t + v \cdot \nabla_x) F_k = B(F_0, F_{k+2}) + B(F_{k+2}, F_0) + \sum_{\substack{i+j = k+2 \,,\\i,j \ge 1}} B(F_i, F_j) \,.
  \end{align}
Together with the well-known $H$-theorem, the order of $\sqrt{\eps}^{-2}$ in \eqref{Interior_Ord_Ana} implies that $F_0 (t,x,v)$ must be a local Maxwellian, hence,
\begin{equation}\label{Local_Maxwellian}
   F_0 (t, x, v) = \M (t, x, v) \,.
\end{equation}
Then the order of $\sqrt{\eps}^0$ in \eqref{Interior_Ord_Ana} implies that $(\rho, \u, T)$ satisfies the compressible Euler equations \eqref{Compressible_Euler_Sys}, and the kinetic part of $F_1$ subjects to \eqref{KinPart-12}. For $k \geq 1$, it is easy to derive from projecting, respectively, the order of $\sqrt{\eps}^k$ in \eqref{Interior_Ord_Ana} into fluid and kinetic parts that the corresponding fluid variables $(\rho_k, u_k, \theta_k)$ obey the linear hyperbolic system
\begin{equation}\label{Linear-Hyperbolic-Syst-k}
  \left\{
    \begin{aligned}
      & \p_t \rho_k + \div_x (\rho u_k + \rho_k \u) = 0 \,, \\
      & \rho \big( \p_t u_k + u_k \cdot \nabla_x \u + \u \cdot \nabla_x u_k \big) - \tfrac{\nabla_x (\rho T)}{\rho} \rho_k + \nabla_x \big( \tfrac{\rho \theta_k + 3 T \rho_k}{3}\big) = \mathcal{F}^\bot_u (\tfrac{F_k}{\sqrt{\M}}) \,, \\
      & \rho \big( \p_t \theta_k + \u \cdot \nabla_x \theta_k + \tfrac{2}{3} \big( \theta_k \div_x \u + 3 T \div_x u_k \big) + 3 u_k \cdot \nabla_x T \big) = \mathcal{F}^\bot_\theta (\tfrac{F_k}{\sqrt{\M}})\,,
    \end{aligned}
  \right.
\end{equation}
and the kinetic parts satisfy
\begin{equation}\label{Ker-Arth-Inte}
(\mathcal{I} - \mathcal{P} ) \big(\tfrac{F_k}{\sqrt{\M}}\big) = \mathcal{L}^{-1} \Big( - \tfrac{(\p_t + v \cdot \nabla_x ) F_{k-2} - \sum_{\substack{i+j = k\,,\\ i, j \ge 1}} B (F_i, F_j) }{\sqrt{\M}} \Big) \,,
\end{equation}
where
\begin{equation*}
  \begin{aligned}
    \mathcal{P} ( \tfrac{F_k}{\sqrt{\M}} ) = \Big\{ \tfrac{\rho_k}{\rho} + u_k \cdot \tfrac{v - \u}{T} + \tfrac{\theta_k}{6 T} ( \tfrac{|v - \u|^2}{T} - 3 ) \Big\} \sqrt{\M} \,,
  \end{aligned}
\end{equation*}
$(\rho, \u, T)$ is the smooth solution to the compressible Euler equations \eqref{Compressible_Euler_Sys}, the source terms $\mathcal{F}^\bot_u (\tfrac{F_k}{\sqrt{\M}})$ and $\mathcal{F}^\bot_\theta (\tfrac{F_k}{\sqrt{\M}})$ are defined as
\begin{equation}
  \begin{aligned}
    & \mathcal{F}^\bot_{u, i} (\tfrac{F_k}{\sqrt{\M}}) = - \sum_{j=1}^3 \p_{x_j} \l T \A_{ij}, \tfrac{F_k}{\sqrt{\M}} \r_v \ (i = 1,2,3) \,, \\
    & \mathcal{F}^\bot_\theta (\tfrac{F_k}{\sqrt{\M}}) = - \div_x \big( 2 T^{\frac{3}{2}} \l \B, \tfrac{F_k}{\sqrt{\M}} \r_v + 2 T \u \cdot \l \A, \tfrac{F_k}{\sqrt{\M}} \r_v \big) - 2 \u \cdot \mathcal{F}^\bot_u (\tfrac{F_k}{\sqrt{\M}})\,.
  \end{aligned}
\end{equation}
Remark that the order of $\sqrt{\eps}^{-1}$ in \eqref{Interior_Ord_Ana} tells us 
\begin{equation*}
	\begin{aligned}
		(\I - \P) ( \tfrac{F_1}{\sqrt{\M}} ) = 0 \,.
	\end{aligned}
\end{equation*}
Before solving the system \eqref{Linear-Hyperbolic-Syst-k}, we will first derive its boundary condition later.

\subsection{Impermeable boundary condition $\u^0 \cdot n = 0$ and diffusive reflection boundary condition}\label{Subsec:2.2}

From substituting the expansion \eqref{Exp1} into the complete diffusive boundary condition \eqref{Diffusive_BC}, the leading order reads
\begin{equation}\label{DFBC-M}
  \begin{aligned}
    \gamma_- \M = M_w \int_{v' \cdot n > 0} (v' \cdot n) \gamma_+ \M \d v' \ \textrm{ on } \Sigma_- \,.
  \end{aligned}
\end{equation}
Once the equation \eqref{DFBC-M}, it is derived from integrating the above equality on $\{ v \cdot n < 0 \}$ with respect to the measure $(v \cdot n) \d v$ that
\begin{equation*}
  \begin{aligned}
    \int_{v \cdot n < 0} (v \cdot n) \M^0 \d v = & \int_{v \cdot n < 0} (v \cdot n) M_w \int_{v^\prime \cdot n >0} (v^\prime \cdot n ) \M^0 \d v^\prime \d v \\
    = & \int_{v \cdot n < 0} (v \cdot n) M_w \d v \int_{v^\prime \cdot n >0} (v^\prime \cdot n ) \M^0 \d v^\prime \\
    = & - \int_{v^\prime \cdot n >0} (v^\prime \cdot n ) \M^0 \d v^\prime \,,
  \end{aligned}
\end{equation*}
where the condition \eqref{Mw-Condition} has been utilized. Consequently,
\begin{equation*}
  \begin{aligned}
    \int_{\R^3} (v \cdot n) \M^0 \d v = 0 \,.
  \end{aligned}
\end{equation*}
Observing $\int_{\R^3} [(v - \u^0) \cdot n] \M^0 \d v = 0$, one therefore has
\begin{equation*}
  \begin{aligned}
    0 = \u^0 \cdot n \int_{\R^3} \M^0 \d v = \rho^0 \u^0 \cdot n \,,
  \end{aligned}
\end{equation*}
hence $\u^0 \cdot n = 0$, which is the boundary condition of the compressible Euler system \eqref{Compressible_Euler_Sys}.

Conversely, under the assumption \eqref{uw-Tw-match}, one has
\begin{equation*}
	\begin{aligned}
		M_w = \tfrac{1}{\rho^0} \sqrt{\tfrac{2 \pi}{T^0}} \M^0 \,.
	\end{aligned}
\end{equation*}
Together with $\u^0 \cdot n = 0$, a direct calculation shows
\begin{equation*}
	\begin{aligned}
		\int_{v' \cdot n > 0} (v' \cdot n) \M^0 \d v' = & \tfrac{\rho^0}{ ( 2 \pi T^0 )^\frac{3}{2} } \int_{v \cdot n > 0} ( v \cdot n ) e^{ - \frac{|v - \u^0|^2}{2 T^0} } \d v \\
		= & \tfrac{\rho^0}{ ( 2 \pi T^0 )^\frac{3}{2} } \int_{( v - \u^0 ) \cdot n > 0} ( ( v - \u^0 ) \cdot n ) e^{ - \frac{|v - \u^0|^2}{2 T^0} } \d v \\
		= & \tfrac{\rho^0}{ ( 2 \pi T^0 )^\frac{3}{2} } \int_{\R^2} \int_{- \infty}^0 (- v_3) e^{- \frac{|\bar{v}|^2 + v_3^2}{2 T^0}} \d \bar{v} \d v_3 \\
		= & \rho^0 \sqrt{\tfrac{T^0}{2 \pi}} \,.
	\end{aligned}
\end{equation*}
As a result, on $\Sigma_-$,
\begin{equation*}
	\begin{aligned}
		M_w \int_{v' \cdot n > 0} (v' \cdot n) \gamma_+ \M \d v' = M_w \int_{v' \cdot n > 0} (v' \cdot n) \M^0 \d v' = M_w \cdot \rho^0 \sqrt{\tfrac{T^0}{2 \pi}} = \M^0 \,,
	\end{aligned}
\end{equation*}
which means that the equation \eqref{DFBC-M} holds under the assumption \eqref{uw-Tw-match}. In other words, the fluid boundary condition $\u^0 \cdot n = 0$ are equivalent to the diffusive reflection boundary condition \eqref{Diffusive_BC} at the leading order, which indicates that there are no strong boundary layers.

\subsection{Viscous boundary layer expansion}\label{Subsec:VBL}

In viscous layer, the scaled normal coordinate is needed:
\begin{equation}
	\begin{aligned}
		\zeta = \tfrac{x_3}{\sqrt{\eps}} \,.
	\end{aligned}
\end{equation}
Dut to no strong boundary layers required, the viscous boundary layer expansion has the form
\begin{equation}
	\begin{aligned}
		F^b_\eps (t, \bar{x}, \zeta) \thicksim \sum_{k \geq 1} \sqrt{\eps}^k F^b_k (t, \bar{x}, \zeta, v) \,,
	\end{aligned}
\end{equation}
where, throughout our paper, the far field condition is always assumed:
\begin{equation}
	\begin{aligned}
		F^b_k (t, \bar{x}, \zeta , v) \to 0 \,, \quad \textrm{as } \zeta \to + \infty \,.
	\end{aligned}
\end{equation}
By similar calculation in \cite{GHW-2020} or \cite{JLT-arXiv-2021}, plugging $F_\eps + F^b_\eps$ into the Boltzmann equation \eqref{BE} gives
\begin{align}\label{Order_Ana}
	\no \sqrt{\eps}^{-1}:\quad & 0= B(\mathfrak{M}^{0}, F^b_1) + B(F^b_1, \mathfrak{M}^{0})\,,\\
	\no \sqrt{\eps}^0:\quad & v_3 \cdot \p_\zeta F^b_1 = \big[ B(\mathfrak{M}^{0}, F^b_2) + B(F^b_2, \mathfrak{M}^{0}) \big] + [ B(F^0_1, F^b_1) + B(F^b_1, F^0_1)] \\
	\no & \qquad \qquad + B(F^b_1, F^b_1) + \zeta \big[ B(\mathfrak{M}^{(1)}, F^b_1) +  B(F^b_1, \mathfrak{M}^{(1)})\big] \,,\\
	\no \sqrt{\eps}: \quad & \p_t F^b_1 + \bar{v} \cdot \nabla_{\bar{x}} F^b_1 + v_3 \cdot \p_\zeta F^b_2 =  \big[ B(\mathfrak{M}^{0}, F^b_3) + B(F^b_3, \mathfrak{M}^{0}) \big] \\
	\no & \qquad \quad + \tfrac{\zeta}{1!} \big[ B(\mathfrak{M}^{(1)}, F^b_2) + B(F^b_2, \mathfrak{M}^{(1)}) \big]  + \tfrac{\zeta^2}{2!} \big[ B(\mathfrak{M}^{(2)}, F^b_1) + B(F^b_1, \mathfrak{M}^{(2)}) \big]\\
	\no & \qquad \quad + \big[ B(F^0_1, F^b_2) + B(F^b_2, F^0_1) \big] + \big[ B(F^0_2, F^b_1) + B(F^b_1, F^0_2) \big] \\
	\no & \qquad \quad + \tfrac{\zeta}{1!} \big[ B(F^{(1)}_1, F^b_1) + B(F^b_1, F^{(1)}_1) \big] + \big[ B(F^b_1, F^b_2) + B(F^b_2, F^b_1) \big]\,,\\
	\no & \cdots \cdots\\
	\no \sqrt{\eps}^k: \quad & \p_t F^b_k + \bar{v} \cdot \nabla_{\bar{x}} F^b_k + v_3 \cdot \p_\zeta F^b_{k+1} \\
	\no & \qquad = \sum_{\substack{i+j=k+2\,,\\ i,j \ge 1 }} B(F^b_i, F^b_j) + \Big[ B(\mathfrak{M}^0, F^b_{k+2}) + B(F^b_{k+2}, \mathfrak{M}^0) \Big]\\
	\no & \qquad \quad + \sum_{\substack{l+j=k+2\,,\\ 1\le l \le N\,, j \ge 1}} \tfrac{\zeta^l}{l!} \Big[ B(\mathfrak{M}^{(l)}, F^b_j) + B(F^b_j, \mathfrak{M}^{(l)}) \Big] \\
	\no & \qquad \quad + \sum_{\substack{i+j=k+2\,,\\ i,j \ge 1 }} \Big[ B(F^0_i, F^b_j) + B(F^b_j, F^0_i) \Big] \\
	& \qquad \quad + \sum_{\substack{i+j+l = k+2\,, \\ 1\le l \le N\,, i, j \ge 1}} \tfrac{\zeta^l}{l!} \Big[ B(F^{(l)}_i, F^b_j) + B(F^b_j, F^{(l)}_i) \Big] \quad {\text{for}}\ k \ge 1 \,,
\end{align}
where the Taylor expansion at $x_3 = 0$ is used:
\begin{equation*}
	\begin{aligned}
		\M = \M^0 + \sum_{1 \leq l \leq N} \tfrac{\zeta^l}{l !} \M^{(l)} + \tfrac{\zeta^{N+1}}{(N+1) !} \widetilde{\M}^{(N+1)} \,, \\
		F_i = F_i^0 + \sum_{1 \leq l \leq N} \tfrac{\zeta^l}{l !} F_i^{(l)} + \tfrac{\zeta^{N+1}}{(N+1) !} \widetilde{F}_i^{(N+1)} \,.
	\end{aligned}
\end{equation*}
Here the number $N \in \mathbb{N}_+$ will be chosen later.

Let
\begin{equation}
	\begin{aligned}
		f^b_k = \tfrac{F^b_k}{\sqrt{\M^0}} \,,
	\end{aligned}
\end{equation}
which can be decomposed as
\begin{equation*}
	\begin{aligned}
		f^b_k = \P^0 f^b_k + (\I - \P^0) f^b_k = \left\{ \tfrac{\rho^b_k}{\rho^0} + u^b_k \cdot \tfrac{v - \u^0}{T^0} + \tfrac{\theta^b_k}{6 T^0} ( \tfrac{|v - \u^0|^2}{T^0} - 3 ) \right\} \sqrt{\M^0} + (\I - \P^0) f^b_k \,.
	\end{aligned}
\end{equation*}
Here $u^b_k = (u^b_{k,1}, u^b_{k,2}, u^b_{k,3}) \in \R^3$. Furthermore, let
\begin{equation}\label{pb_k}
	\begin{aligned}
		p^b_k = \tfrac{\rho^0 \theta^b_k + 3 T^0 \rho^b_k}{3} \,.
	\end{aligned}
\end{equation}
Following the Theorem 1.1 of \cite{GHW-2020}, the following relations are obyed
\begin{equation}\label{ub_13-pb_1}
	\begin{aligned}
		u^b_{1,3} (t, \bar{x}, \zeta) \equiv 0 \,, \ p^b_1 (t, \bar{x}, \zeta) \equiv 0 \,, \ \forall (t, \bar{x}, \zeta) \in [0, \tau] \times \R^2 \times \R_+ \,,
	\end{aligned}
\end{equation}
and $(u^b_{k,1}, u^b_{k,2}, \theta^b_k)$ $(k \geq 1)$ satisfy the following {\em linear compressible Prandtl-type} equations
\begin{equation}\label{Linear-Prandtl}
	\left\{
	\begin{aligned}
		\rho^0 (\p_t + \bar{\u}^0 \cdot \nabla_{\bar{x}}) u^b_{k,i} + \rho^0 (\p_{x_3} \u^0_3 \zeta + u^0_{1,3}) \p_{\zeta} u^b_{k,i} + & \rho^0 \bar{u}^b_k \cdot \nabla_{\bar{x}} \u^0_i + \tfrac{\p_{x_3} p^0}{3 T^0} \theta^b_k \\
		= & \mu (T^0) \p_{\zeta}^2 u^b_{k,i} + \mathsf{f}^b_{k-1, i} \ (i=1,2) \,, \\
		\rho^0 \p_t \theta^b_k + \rho^0 \bar{\u}^0 \cdot \nabla_{\bar{x}} \theta^b_k + \rho^0 \big( \p_{x_3} \u^0 \zeta + \u^0_{1,3} \big) \p_\zeta \theta^b_k + & \tfrac{2}{3} \rho^0 \div_{x} \u^0 \theta^b_k \\
		= & \tfrac{3}{5} \kappa(T^0) \p_{\zeta\zeta} \theta^b_k + \mathsf{g}^b_{k-1} \,, \\
		\lim_{\zeta \to \infty} (\bar{u}^b_k, \theta^b_k) (t, \bar{x}, \zeta) = & 0 \,,
	\end{aligned}
	\right.
\end{equation}
and $(\I - \P^0) f^b_{k+1}$, $u^b_{k+1, 3}$, $p^b_{k+1}$ are determined by the equations
\begin{equation}\label{f_b_k+1-Kernal-Ortho}
	\begin{aligned}
		(\I - \P^0) & f^b_{k+1} = (\mathcal{L}^0)^{-1} \Big\{ - (\I - \P^0) (v_3 \p_\zeta \P^0 f^b_k) + \tfrac{\zeta}{\sqrt{\M^0}} \big[ B(\M^{(1)}, \sqrt{\M^0} \P^0 f^b_k) \\
		& + B(\sqrt{\M^0} \P^0 f^b_k, \M^{(1)})\big] + \tfrac{1}{\sqrt{\M^0}} \big[ B(F^0_1, \sqrt{\M^0} \P^0 f^b_k) + B(\sqrt{\M^0} \P^0 f^b_k, F^0_1) \big] \\
		& + \tfrac{1}{\sqrt{\M^0}} \big[ B(\sqrt{\M^0} f^b_1, \sqrt{\M^0} \P^0 f^b_k) + B(\sqrt{\M^0} \P^0 f^b_k, \sqrt{\M^0} f^b_1) \big] \Big\} + J^b_{k-1}\,,
	\end{aligned}
\end{equation}
and
\begin{equation}\label{u_b_k+1-derivative}
	\begin{aligned}
		\p_\zeta u^b_{k+1, 3} = - \frac{1}{\rho^0} ( \p_t \rho^b_k + \div_{\bar{x}} ( \rho^0 \bar{u}^b_k + \rho^b_k \bar{\u}^0 ) )\,, \ \lim_{\zeta \to \infty} u^b_{k+1,3} (t, \bar{x}, \zeta) = 0 \,,
	\end{aligned}
\end{equation}
and
\begin{equation}\label{p_b_k+1-derivative}
	\left\{
	\begin{aligned}
		\p_\zeta p^b_{k+1}
		= - \rho^0 \p_t u^b_{k,3} - \rho^0 \bar{\u}^0 \!\cdot\! \nabla_{\bar{x}} u^b_{k,3} + \rho^0 \p_{x_3} \u^0_3 u^b_{k,3} + \tfrac{4}{3} \mu(T^0) \p_{\zeta\zeta} u^b_{k,3} \\
		- \tfrac{4}{3} \rho^0 \p_\zeta \Big[ \big(\p_{x_3} \u^0 \zeta + u^0_{1,3} \big) u^b_{k,3} \Big] - \p_\zeta \l T^0 \A^0_{33}, J^b_{k-1} \r + W^b_{k-1,3} \,, \\
		\lim_{\zeta \to \infty} \theta^b_{k+1} (t, \bar{x}, \zeta) = 0 \,,
	\end{aligned}
	\right.
\end{equation}
where the source terms $\mathsf{f}^b_{k-1, i} \ (i = 1, 2)$ and $\mathsf{g}^b_{k-1}$ are
\begin{equation}\label{fg-k}
	\begin{aligned}
		\mathsf{f}^b_{k-1, i} = -\rho^0 \p_\zeta [ ( \p_{x_3} \u^0_i \zeta + u^0_{1,i} + u^b_{1,i} ) u^b_{k,3} ] - ( \p_{x_i} - \tfrac{\p_{x_i} p^0}{p^0} ) p^b_k \\
		+ W^b_{k-1,i} - T^0 \p_\zeta \l J^b_{k-1}, \A^0_{3i} \r \,, \\
		\mathsf{g}^b_{k-1} = - \rho^0 \p_\zeta \Big[ \big( 3 \p_\zeta T^0 \zeta + \theta^0_1 + \theta^b_1 \big) u^b_{k,3}\Big] + \tfrac{3}{5} H^b_{k-1} - \tfrac{6}{5} (T^0)^{\frac{3}{2}} \p_\zeta \l J^b_{k-1}, \B^0_3 \r\\
		+ \tfrac{3}{5} \big\{ 2 \p_t + 2 \bar{\u} ^0 \cdot \nabla_{\bar{x}} + \tfrac{10}{3} \div_{x} \u^0 \big\} p^b_k \,,
	\end{aligned}
\end{equation}
and
\begin{align}\label{W_b_k-1}
	W^b_{k-1,i} = - \sum_{j=1}^2 \p_{x_j} \l T^0 (\I- \P^0) f^b_k, \A^0_{ij} \r \,, \qquad {\text{for}}\ i =1,2,3\,,
\end{align}
\begin{align}\label{H_b_k-1}
	H^b_{k-1} = - \sum_{j=1}^2 \p_{x_j} \l (T^0)^{\frac{3}{2}} \B^0_j + \sum_{l=1}^2 2 T^0 \u^0_l \A^0_{jl},  (\mathcal{I} - \P^0) f^b_k \r - 2 \bar{\u}^0 \cdot \bar{W}^b_{k-1}\,,
\end{align}
and
\begin{align}\label{J_b_k-1}
	\no J^b_{k-1} = & (\mathcal{L}^0)^{-1} \Big\{ - (\mathcal{I}-\mathcal{P}^0) \Big( \tfrac{1}{\sqrt{\M^0}} \big\{ \p_t + \bar{v} \cdot \nabla_{\bar{x}} \big\} F^b_{k-1} \Big) - (\mathcal{I}-\mathcal{P}^0) \big( v_3 \p_\zeta ( \mathcal{I} - \mathcal{P}^0 \big) f^b_k \big) \\
	\no & + \sum_{\substack{l+j=k+1\,,\\2 \le l \le N, j\ge 1}} \tfrac{\zeta^l}{l!} \tfrac{1}{\sqrt{\M^0}} \big[ B(\M^{(l)}, \sqrt{\M^0} f^b_j) + B(\sqrt{\M^0} f^b_j, \M^{(l)}) \big] \\
	\no & + \sum_{\substack{l+j=k+1\,,\\i\ge 2, j\ge 1}} \tfrac{1}{\sqrt{\M^0}} \big[ B(F^0_i, \sqrt{\M^0} f^b_j) + B(\sqrt{\M^0} f^b_j, F^0_i) \big] \\
	\no & + \sum_{\substack{l+j=k+1\,,\\ i,j\ge 2}} \tfrac{1}{\sqrt{\M^0}} \big[ B(\sqrt{\M^0} f^b_i, \sqrt{\M^0} f^b_j) + B(\sqrt{\M^0} f^b_j, \sqrt{\M^0} f^b_i) \big] \\
	\no & + \sum_{\substack{i+j+k = k+1\,,\\1\le l \le N, i,j \ge 1}} \tfrac{1}{\sqrt{\M^0}} \tfrac{\zeta^l}{l!} \big[ B(F^{(l)}_i, \sqrt{\M^0} f^b_j) + B(\sqrt{\M^0} f^b_j, F^{(l)}_i)\big] \\
	\no & + \tfrac{\zeta}{\sqrt{\M^0}} \big[ B(\M^{(1)}, \sqrt{\M^0} (\mathcal{I} - \mathcal{P}^0) f^b_k) + B(\sqrt{\M^0} (\mathcal{I} - \mathcal{P}^0) f^b_k, \M^{(1)}) \big] \\
	\no & + \tfrac{1}{\sqrt{\M^0}} \big[ B(F^0_1, \sqrt{\M^0} (\mathcal{I} - \mathcal{P}^0) f^b_k) + B(\sqrt{\M^0} (\mathcal{I} - \mathcal{P}^0) f^b_k, F^0_1) \big] \\
	& + \tfrac{1}{\sqrt{\M^0}} \big[ B(\sqrt{\M^0} f^b_1, \sqrt{\M^0} (\mathcal{I} - \mathcal{P}^0) f^b_k) + B(\sqrt{\M^0} (\mathcal{I} - \mathcal{P}^0) f^b_k, \sqrt{\M^0} f^b_1) \big] \Big\} \,.
\end{align}
We emphasize that $W^b_{k-1}$, $H^b_{k-1}$ and $J^b_{k-1}$ depend on $f^b_j$ $(1 \le j \le k-1)$. Moreover, when $k=1$,  $J^b_0 = W^b_0 = H^b_0 = \mathsf{f}^b_0 = \mathsf{g}^b_0 = 0$. 

Here the initial conditions of \eqref{Linear-Prandtl} are imposed on
\begin{equation}\label{IC-Prandtl}
	\begin{aligned}
		(\bar{u}^b_k, \theta^b_k) (0, \bar{x}, \zeta) = (\bar{u}^{b, in}_k , \theta^{b, in}_k ) (\bar{x}, \zeta) \in \R^2 \times \R \,, \quad k = 1,2,3 \,, \cdots
	\end{aligned}
\end{equation}
with the compatibility $\lim_{\zeta \to \infty} (\bar{u}^{b, in}_k , \theta^{b, in}_k) (\bar{x}, \zeta) = 0$. Before solving the linear Prandtl system \eqref{Linear-Prandtl}, its boundary conditions should be initially derived later.

\subsection{Knudsen boundary layers expansion}

In the Knudsen layer, the new scaled normal coordinate is introduced:
\begin{equation*}
	\begin{aligned}
		\xi = \tfrac{x_3}{\epsilon} \,.
	\end{aligned}
\end{equation*}
Together with \eqref{Interior_Ord_Ana}, it is derived from plugging the expansion \eqref{Exp1} into \eqref{BE} that
\begin{equation}\label{KL_Order_Ana}
  \begin{aligned}
   \sqrt{\eps}^{-1} :\qquad & v_3 \cdot \p_\xi F^{bb}_1 = B(\M^0, F^{bb}_1) + B(F^{bb}_1, \M^0)\\
   \sqrt{\eps}^0 : \qquad & v_3 \cdot \p_\xi F^{bb}_2 - \big[ B(\M^0, F^{bb}_2) + B(F^{bb}_2, \M^0) \big] \\  = & B(F^0_1 + F^{b,0}_1, F^{bb}_1) + B(F^{bb}_1, F^0_1 + F^{b,0}_1) + B(F^{bb}_1, F^{bb}_1)\\
   ......&\\
   \sqrt{\eps}^k :\qquad & v_3 \cdot \p_\xi F^{bb}_{k+2} - \big[ B(\M^0, F^{bb}_{k+2}) + B(F^{bb}_{k+2}, \M^0) \big] \\
   = & - \big\{ \p_t + \bar{v} \cdot \nabla_{\bar{x}}
   \big\} F^{bb}_k + \sum_{\substack{j+2l=k+2\,,\\ 1 \le l \le N, j\ge 1}} \frac{\xi^l}{l!} \big[ B(\M^{(l)}, F^{bb}_j) + B(F^{bb}_j, \M^{(l)})\big] \\
   & + \sum_{\substack{i+j=k+2\,,\\ i, j\ge 1}} \big[ B(F^0_i + F^{b,0}_i, F^{bb}_j) + B(F^{bb}_j, F^0_i + F^{b,0}_i) + B(F^{bb}_i, F^{bb}_j) \big] \\
   & + \sum_{\substack{i+2l+j=k+2\,,\\1\le l \le N, i,j\ge 1}} \frac{\xi^l}{l!} \big[ B(F^{(l)}_i, F^{bb}_j) + B(F^{bb}_j, F^{(l)}_i)\big] \\
   & + \sum_{\substack{i+l+j=k+2\,,\\1\le l \le N, i,j\ge 1}} \frac{\xi^l}{l!} \big[ B(F^{b, (l)}_i, F^{bb}_j) + B(F^{bb}_j, F^{b, (l)}_i)\big]\,.
  \end{aligned}
\end{equation}
Here the Taylor expansions of $\M$, $F_k$ at $x_3=0$ and $F^b_i$ at $\zeta = 0$ are used:
\begin{equation*}
  \begin{aligned}
   \M = & \M^0 + \sum_{1 \le l \le N} \tfrac{\xi^l}{l !} \M^{(l)} + \tfrac{\xi^{N+1}}{(N+1)!} \widetilde{\M}^{(N+1)}\,,\\
   F_i = & F_i^0 + \sum_{1 \le l \le N} \tfrac{\xi^l}{l !} F_i^{(l)} + \tfrac{\xi^{N+1}}{(N+1)!} \widetilde{F}_i^{(N+1)}\,, \\
   F_i^b = & F_i^{b,0} + \sum_{1 \leq l \leq N} \tfrac{\xi^l}{l!} F_i^{b, (l)} + \tfrac{\xi^{N+1}}{(N+1)!} \widetilde{F}_i^{b, (N+1)} \,,
  \end{aligned}
\end{equation*}
where $N=3$ in the truncated expansion form \eqref{Hilbert_Expansion1}. 

Let $f_k = \frac{F_k}{\sqrt{\M}}$ and $f^{bb}_k = \frac{F^{bb}_k}{\sqrt{\M^0}}$. \eqref{KL_Order_Ana} implies
\begin{equation}\label{Kudesen_Layer_Eq}
 \left\{
  \begin{aligned}
   &v_3 \p_{\xi} f^{bb}_k + \mathcal{L}^0 f^{bb}_k = \mathcal{S}^{bb}_k \,,\\
   &\lim_{\xi \to +\infty} f^{bb}_k (t, \bar{x}, \xi, v) = 0\,,
  \end{aligned}
  \right.
\end{equation}
where $S^{bb}_k = S^{bb}_{k,1} + S^{bb}_{k,2}$ with
 \begin{align}\label{S_bb_k+1_1+2}
   \no S^{bb}_{k,1} =& - \mathcal{P}^0 \Big\{ \frac{(\p_t + \bar{v} \cdot \nabla_{\bar{x}}) F^{bb}_{k-2}}{\sqrt{\M^0}} \Big\} \in \mathcal{N}^0 \,, \\
   \no S^{bb}_{k,2} = & \sum_{\substack{j+2l=k\,,\\ 1 \le l \le N, j\ge 1}} \frac{\xi^l}{l!} \frac{1}{\sqrt{\M^0}} \big[ B(\M^{(l)}, \sqrt{\M^0} f^{bb}_j) + B(\sqrt{\M^0} f^{bb}_j, \M^{(l)})\big] \\
   \no & + \sum_{\substack{i+j=k\,,\\ i, j\ge 1}} \frac{1}{\sqrt{\M^0}} \big[ B(F^0_i + F^{b,0}_i, \sqrt{\M^0} f^{bb}_j) + B(\sqrt{\M^0} f^{bb}_j, F^0_i + F^{b,0}_i)\big] \\
   & + \sum_{\substack{i+2l+j=k\,,\\1\le l \le N, i,j\ge 1}} \frac{\xi^l}{l!} \frac{1}{\sqrt{\M^0}} \big[ B(F^{(l)}_i, \sqrt{\M^0} f^{bb}_j) + B(\sqrt{\M^0} f^{bb}_j, F^{(l)}_i)\big] \\
   \no & + \sum_{\substack{i+l+j=k\,,\\1\le l \le N, i,j\ge 1}} \frac{\xi^l}{l!} \frac{1}{\sqrt{\M^0}} \big[ B(F^{b, (l)}_i, \sqrt{\M^0} f^{bb}_j) + B(\sqrt{\M^0} f^{bb}_j, F^{b, (l)}_i)\big] \\
   \no & + \sum_{\substack{i+j=k\,,\\ i,j\ge 1}} \frac{1}{\sqrt{\M^0}} B(\sqrt{\M^0} f^{bb}_i, \sqrt{\M^0} f^{bb}_j) - (\mathcal{I} - \mathcal{P}^0) \Big\{ \frac{\{\p_t + \bar{v} \cdot \nabla_{\bar{x}} \} F^{bb}_{k-2}}{\sqrt{\M^0}} \Big\} \in (\mathcal{N}^0)^\bot \,.
 \end{align}
Here  the notation $F^{bb}_{-1} = F^{bb}_0=0 $ are used, and
\begin{align*}
	S^{bb}_1 = S^{bb}_{1,1} = S^{bb}_{1,2} =0\,,\ S^{bb}_{2,1} = \mathcal{P}^0 S^{bb}_2 =0 \,.
\end{align*}

Now, the following lemma is quoted to deal with the source term $S^{bb}_{k,1} \in \mathcal{N}^0$.

\begin{lemma}[\cite{Bardos-Caflisch-Nicolaenko-1986-CPAM}]\label{Lmm-fbb-k1}
	Assume that
	\begin{equation}\label{S_bb_k1}
	\begin{aligned}
	S^{bb}_{k,1}= \big\{ a_k + b_k \cdot (v- \u^0) + c_k |v-\u^0|^2 \big\} \sqrt{\M^0}
	\end{aligned}
	\end{equation}
	satisfies
	\begin{align*}
	\lim_{\xi \to \infty} e^{\eta \xi} |(a_k, b_k, c_k) (t, \bar{x}, \xi)| = 0
	\end{align*}
	for some positive constant $\eta > 0$. Then there exists a function
	\begin{equation*}
	f^{bb}_{k,1} = \big\{ \Psi_k v_3 + \Phi_{k,1} v_3 (v_1 - \u^0_1) + \Phi_{k,2} v_3 (v_2 - \u^0_2) + \Phi_{k,3} + \Theta_k v_3 |v- \u^0|^2 \big\} \sqrt{\M^0} 	
	\end{equation*}
	such that
	\begin{equation*}
	\begin{aligned}
	v_3 \p_\xi f^{bb}_{k,1} - S^{bb}_{k,1} \in (\mathcal{N}^0)^\bot \,,
	\end{aligned}
	\end{equation*}
	where
	\begin{equation}\label{f_bb_k1_Coef}
	\begin{array}{l}
	\Psi_k (t, \bar{x}, \xi)
	= - \int_\xi^{+\infty} \big( \tfrac{2}{T^0} a_k + 3 c_k \big) (t, \bar{x}, s) \d s\,, \\ [1.5mm]
	\Phi_{k,i} (t, \bar{x}, \xi)
	= - \int_\xi^{+\infty} \tfrac{1}{T^0} b_{k,i} (t, \bar{x}, s) \d s\,,\ i=1,2\,, \\ [1.5mm]
	\Phi_{k,3} (t, \bar{x}, \xi)
	= - \int_\xi^{+\infty} b_{k,3} (t, \bar{x}, s) \d s\,, \\ [1.5mm]
	\Theta_k (t, \bar{x}, \xi) = \tfrac{1}{5(T^0)^2} \int_\xi^{+\infty} a_k (t, \bar{x}, s) \d s\,.
	\end{array}
	\end{equation}
	Moreover, there holds
	\begin{equation*}
	\begin{aligned}
	& |v_3 \p_\xi f^{bb}_{k,1} - S^{bb}_{k,1}| \leq C |(a_k, b_k, c_k) (t, \bar{x}, \xi)| \l v \r^4 \sqrt{\M^0} \,, \\
	& |f^{bb}_{k,1} (t, \bar{x}, \xi, v)| \leq C \l v \r^3 \sqrt{\M^0} \int_\xi^\infty |(a_k, b_k, c_k)| \to 0 \textrm{ as } \xi \to \infty \,.
	\end{aligned}
	\end{equation*}
\end{lemma}

Let $f_{k, 2}^{bb} = f^{bb}_k - f^{bb}_{k, 1}$. One thereby obtains the following linear Knudsen boundary layer equation 
\begin{equation}\label{Knudsen-Layer-1}
\left\{
\begin{aligned}
& v_3 \p_\xi f^{bb}_{k, 2} + \mathcal{L}^0 f^{bb}_{k, 2} = S^{bb}_{k, 2} - \big( v_3 \p_\xi f^{bb}_{k,1} - S^{bb}_{k, 1} + \mathcal{L}^0 f^{bb}_{k,1} \big) \in (\mathcal{N}^0 )^\perp \,, \\
& \lim_{\xi \to +\infty} f^{bb}_{k,2} (t, \bar{x}, \xi, v) =0 \,.
\end{aligned}
\right.
\end{equation}
Before solving the problem \eqref{Knudsen-Layer-1}, its boundary condition on $\Sigma_-$ should be first derived later.

\subsection{Constructions of boudary conditions for all orders}

By plugging \eqref{Exp1} into the complete diffusive boundary condition \eqref{Diffusive_BC},
\begin{equation}\label{Diffusive_BC_Order_Ana}
	\begin{aligned}
		\gamma_- (F_k + F^b_k + F^{bb}_k) = M_w \int_{v^\prime \cdot n >0} (v^\prime \cdot n) \gamma_+ (F_k + F^b_k + F^{bb}_k) \d v^\prime \quad {\text{on }} \Sigma_- \,.
	\end{aligned}
\end{equation}
for $k \geq 1$. We now construct the boundary conditions for the linear hyperbolic system \eqref{Linear-Hyperbolic-Syst-k}, the linear Prandtl system \eqref{Linear-Prandtl} and the Knudsen boundary layer equation \eqref{Knudsen-Layer-1} such that the complete diffusive boundary conditions \eqref{Diffusive_BC_Order_Ana} holds for all $k \geq 1$.

We first rewrite \eqref{Diffusive_BC_Order_Ana} as
\begin{equation}\label{KL-BC}
	\begin{aligned}
		f_{k,2}^{bb} (t, \bar{x}, 0, v) |_{v_3 > 0} = \mathcal{D}_w f_{k,2}^{bb} (t, \bar{x}, 0, v) + \mathbbm{f}_k (t, \bar{x}, v) \,,
	\end{aligned}
\end{equation}
where the operator $\mathcal{D}_w$ is given in \eqref{Dw}, and $\mathbbm{f}_k (t, \bar{x}, v)$ is defined in \eqref{fk}, namely,
\begin{equation*}
	\begin{aligned}
		\mathbbm{f}_k (t, \bar{x}, v) = \mathbbm{f}_k (t, \bar{x}, v) = - ( \I - \mathcal{D}_w ) ( f_k + f^b_k + f^{bb}_{k,1} ) (t, \bar{x}, 0, v) \,.
	\end{aligned}
\end{equation*}

From Coron, Golse and Sulem's work \cite{Coron-Golse-Sulem-1988-CPAM}, the solvability of the Knudsen boundary layer equation \eqref{Knudsen-Layer-1} with boundary condition \eqref{KL-BC} is
\begin{equation}\label{KL_Sol_Condition}
  \begin{aligned}
    \int_{v \cdot n < 0} ( v \cdot n ) \mathbbm{f}_k (t, \bar{x}, v) \sqrt{\M^0} \d v = 0 \,.
  \end{aligned}
\end{equation}
Recalling the condition \eqref{Mw-Condition}, i.e., $\int_{v \cdot n < 0} |v \cdot n| M_w (t,x,v) \d v = 1$, one derives from \eqref{KL_Sol_Condition} that
\begin{equation*}
  \begin{aligned}
    \int_{v \cdot n < 0} (v \cdot n) ( f_k + f^b_k + f^{bb}_{k,1} ) \sqrt{\M^0} \d v = \int_{v \cdot n < 0} (v \cdot n) \mathcal{D}_w ( f_k + f^b_k + f^{bb}_{k,1} ) \sqrt{\M^0} \d v \\
    = \int_{v \cdot n < 0} (v \cdot n) \tfrac{M_w}{\sqrt{\M^0}} \int_{v' \cdot n > 0} (v' \cdot n) ( f_k + f^b_k + f^{bb}_{k,1} ) \sqrt{\M^0} \d v' \sqrt{\M^0} \d v \\
    = - \int_{v' \cdot n > 0} (v' \cdot n) ( f_k + f^b_k + f^{bb}_{k,1} ) \sqrt{\M^0} \d v' \,,
  \end{aligned}
\end{equation*}
which means that the solvability \eqref{KL_Sol_Condition} is equivalent to
\begin{equation}\label{Solva_Condition_f_bb_k_2}
  \begin{aligned}
    \int_{\R^3} (v \cdot n) ( f_k + f^b_k + f^{bb}_{k,1} ) (t, \bar{x}, 0 ,v) \sqrt{\M^0} \d v = 0 \,.
  \end{aligned}
\end{equation}

{\bf Claim:} {\em The solvability \eqref{Solva_Condition_f_bb_k_2} implies the boundary condition of the linear hyperbolic system \eqref{Linear-Hyperbolic-Syst}:
\begin{equation}\label{High_Ord_Fluid_BC}
\begin{aligned}
  u^0_k \cdot n = & T^0 (\Psi_k + 5 T^0 \Theta_k) (t, \bar{x}, 0) - \tfrac{1}{\rho^0} \int_{\R^3} (v \cdot n) (\mathcal{I} - \mathcal{P}^0) f_k (t, \bar{x}, 0, v) \sqrt{\M^0} \d v \\
  & - \tfrac{1}{\rho^0} \int_0^\infty ( \partial_t \rho^b_{k-1} + \div_{\bar{x}} ( \rho^0 \bar{u}_{k-1}^b + \rho^b_{k-1} \bar{\u}^0 ) ) (t, \bar{x}, \zeta) \d \zeta \,.
\end{aligned}
\end{equation}
}

Indeed, after splitting $f_k$ and $f^b_k$ as
\begin{equation*}
  \begin{aligned}
   f_k = & \big\{ \tfrac{\rho_k}{\rho} + u_k \cdot \tfrac{v - \u}{T} + \tfrac{3 \theta_k}{6 T} \big( \tfrac{|v-\u|^2}{T} -3 \big) \big\} \sqrt{\M} + (\mathcal{I} - \mathcal{P}) f_k \,, \\
   f^b_k = & \big\{ \tfrac{\rho_k^b}{\rho} + u_k^b \cdot \tfrac{v - \u^0}{T} + \tfrac{3 \theta_k^b}{6 T} \big( \tfrac{|v-\u^0|^2}{T} -3 \big) \big\} \sqrt{\M} + (\mathcal{I} - \mathcal{P}^0) f_k^b \,,
  \end{aligned}
\end{equation*}
the integral $\int_{\R^3} (v \cdot n) f_k (t, \bar{x}, 0, v) \sqrt{\M^0} \d v$ can be calculated by
\begin{align*}
  & \int_{\R^3} (v \cdot n) ( f_k + f^b_k ) (t, \bar{x}, 0, v) \sqrt{\M^0} \d v \\
  = & \int_{\R^3} (v \cdot n) (\mathcal{I} - \mathcal{P}^0) ( f_k + f^b_k ) (t, \bar{x}, 0, v) \sqrt{\M^0} \d v \\
  & + \int_{\R^3} (v - \u^0 ) \cdot n \big\{ \tfrac{\rho_k^0 + \rho^{b, 0}_k}{\rho^0} + ( u_k^0 + u_k^{b,0} ) \cdot \tfrac{v - \u^0}{T^0} + \tfrac{3 ( \theta_k^0 + \theta_k^{b,0} ) }{6 T^0} \big( \tfrac{|v-\u^0|^2}{T^0} -3 \big) \big\} \M^0 \d v \\
  = & \rho^0 ( u^0_k + u_k^{b,0} ) \cdot n + \int_{\R^3} (v \cdot n) (\mathcal{I} - \mathcal{P}^0) ( f_k + f^b_k ) (t, \bar{x}, 0, v) \sqrt{\M^0} \d v \,,
\end{align*}
where $\u^0 \cdot n =0$ has been used. On the other hand, by $\u^0 \cdot n = 0$ and the definition of $f^{bb}_{k,1}$ in Lemma \ref{Lmm-fbb-k1}, there holds
\begin{align*}
   & \int_{\R^3} (v \cdot n) f^{bb}_{k, 1} (t, \bar{x}, 0, v) \sqrt{\M^0} \d v \\
   = & \int_{\R^3} [(v - \u^0) \cdot n] \big( \Psi_k v_3 + \sum_{i=1,2} \Phi_{k, i} v_3 (v_i - \u^0_i) + \Phi_{k, 3} + \Theta_k v_3 |v-\u^0|^2 \big) \M^0 \d v \\
   = & - \rho^0 T^0 (\Psi_k + 5 T^0 \Theta_k)\,.
\end{align*}
Collecting the above equations reduces 
\begin{equation*}
	\begin{aligned}
		( u^0_k + u_k^{b,0} ) \cdot n = T^0 (\Psi_k + 5 T^0 \Theta_k) (t, \bar{x}, 0) - \tfrac{1}{\rho^0} \int_{\R^3} (v \cdot n) (\mathcal{I} - \mathcal{P}^0) f_k (t, \bar{x}, 0, v) \sqrt{\M^0} \d v \,.
	\end{aligned}
\end{equation*}
By \eqref{u_b_k+1-derivative},
\begin{equation*}
	\begin{aligned}
		u_k^{b,0} \cdot n = \tfrac{1}{\rho^0} \int_0^\infty ( \partial_t \rho^b_{k-1} + \div_{\bar{x}} ( \rho^0 \bar{u}_{k-1}^b + \rho^b_{k-1} \bar{\u}^0 ) ) (t, \bar{x}, \zeta) \d \zeta \,.
	\end{aligned}
\end{equation*}
Consequently, the claim \eqref{High_Ord_Fluid_BC} holds.

By now, the claim \eqref{High_Ord_Fluid_BC} ensures that $(\rho_k, u_k, \theta_k)$ is solved, which means that $f_k$ is completely determined. There is only $f^b_k$ undetermined in the expression of $\mathbbm{f}_k (t, \bar{x}, v)$. We now impose the following boundary values on the linear Prandtl system \eqref{Linear-Prandtl}:
\begin{equation}\label{BC-Prandtl}
	\begin{aligned}
		\bar{u}^b_k (t, \bar{x}, 0) = \bar{u}_w (t, \bar{x}) \,, \quad \theta^b_k (t, \bar{x}, 0) = T_w (t, \bar{x}) \,.
	\end{aligned}
\end{equation}
The functions $(\bar{u}^b_k, \theta^b_k) (t, \bar{x}, \zeta)$ is thereby solved by \eqref{Linear-Prandtl}, which indicates that $f_k^b$ is completely determined. As a result, the function $\mathbbm{f}_k (t, \bar{x}, v)$ is completely determined. At the end, we imposed the boundary condition \eqref{KL-BC} on the problem \eqref{Knudsen-Layer-1}. Then $f^{bb}_{k,2} (t, \bar{x}, \xi, v)$ can be solved, and the equation \eqref{Diffusive_BC_Order_Ana} holds for $k \geq 1$.

\subsection{Truncation for Hilbert expansion}

In this subsection, the truncation of the Hilbert expansion \eqref{Exp1} will be done. Our goal is to prove the compressible Euler limit from the scaled Boltzmann equation by above Hilbert type expansion. The key point is to prove that the remainders of expansion will go to zero as the Knudsen number $\eps \to 0$. Mathematically, this means we search for a special class of solutions of the original scaled Boltzmann equation for sufficiently small Knudsen number $\eps$. Since the more terms are expanded, the more special the solutions are. We hope that the terms in the expansion are as less as possible.

Our truncation form is taken as
\begin{equation}\label{Exp-truncated}
	\begin{aligned}
		F_\eps (t, x, v) = \M (t,x,v) + \sum_{k = 1}^6 \sqrt{\eps}^k \big[ F_k (t,x,v) + F^b_k (t, \bar{x}, \tfrac{x_3}{\sqrt{\eps}}, v) + F^{bb}_k (t, \bar{x}, \tfrac{x_3}{\eps}, v) \big] \\
		+ \sqrt{\eps}^5 F_{R, \eps} (t,x,v) \,,
	\end{aligned}
\end{equation}
where $F_{R,\eps}$ is the remainder term.

First, $\M$ in \eqref{Max-M} is governed by the compressible Euler system \eqref{Compressible_Euler_Sys}-\eqref{BC-CEuler}-\eqref{IC-CEuler}. 

Second, the interiors $ F_1 = \P ( \frac{F_1}{\sqrt{\M}} ) \sqrt{\M} $ and $ F_k = \{ \P ( \frac{F_k}{\sqrt{\M}} ) + ( \I - \P ) ( \frac{F_k}{\sqrt{\M}} ) \} \sqrt{\M} $ ($k = 2,3,4$) are determined as follows: $( \I - \P ) ( \frac{F_k}{\sqrt{\M}} )$ ($k = 2,3,4$) are given in \eqref{Ker-Arth-Inte}, and the fluid variables $(\rho_k, u_k, \theta_k)$ ($k =1, 2,3,4$) associated with $\P ( \frac{F_k}{\sqrt{\M}} )$ are the solutions to the linear hyperbolic system \eqref{Linear-Hyperbolic-Syst}. Then $F_k$ ($k = 1,2,3,4$) are completed determined. 

Third, the viscous layer expansions $F^b_1 = (\P^0 f_1^b ) \sqrt{\M^0} $ and $ F^b_k = \{ \P^0 f^b_k + ( \I - \P^0 ) f^b_k \} \sqrt{\M^0} $ $(k=2,3,4)$ are determined as follows: The kinetic part $ ( \I - \P^0 ) f^b_k $ $(k=2,3,4)$ are represented by \eqref{f_b_k+1-Kernal-Ortho}, and the fluid variables $(\rho^b_k, u^b_k, \theta^b_k)$ corresponding to $\P^0 f^b_k$ $(k=1, 2,3,4)$ are determined by the equations \eqref{u_b_k+1-derivative}-\eqref{p_b_k+1-derivative} and the linear Prandtl system \eqref{Linear-Prandtl}. For $k = 5,6$, we choose $F^b_k$ satisfying their fluid variables vanished, i.e., $F^b_k = [ (\I - \P^0) f^b_k ] \sqrt{\M^0}$, which are expressed by \eqref{f_b_k+1-Kernal-Ortho}, which mean that they are not completely determined.

Forth, the Knudsen boundary layers $F^{bb}_k$ $(k = 1,2,3,4)$ are determined by \eqref{Kudesen_Layer_Eq} and \eqref{KL-BC}, who depends the expansions $F^{bb}_i$ ($i \leq k-1$) and $F_j, F^b_j$ $(j \leq k)$. As stated above, we know that $F^{bb}_k$ $(k = 1,2,3,4)$ are completely determined. For $k=5,6$, $F^{bb}_k$ can be solved by the similar way of $F^{bb}_i$ ($i \leq 4$). Since $F_j, F^b_j$ $(j=5,6)$ are only partially determined, $F^{bb}_k$ $(k=5,6)$ are also not completely determined. We also emphasize that the number $N \in \mathbb{N}_+$ appeared in the Taylor expansions before will be chosen by $N = 5$ for the truncated Hilbert expansion \eqref{Exp-truncated}.

At the end, we obtain the  remainder $F_{R,\eps}$ with order $\sqrt{\eps}^5$, who obeys the equations \eqref{Remainder_Eq} with the complete diffusive boundary condition \eqref{Remainder_BC}.

\section{Uniform bounds for expanded coefficients}\label{Sec:F1F2}

In this section, the uniform bounds for the expansion terms $F_k$, $F^b_k$ and $F^{bb}_k$ ($1 \leq k \leq 6$) shall be derived, which will be applied to control the coefficients in the remainder equation \eqref{Remainder_Eq}. 

First, as shown in Section \ref{Sec:FormalAnaly} above, the fluid variables $(\rho_k, u_k, \theta_k)$ of $F_k$ ($k=1,2,3,4$) subject to the linear hyperbolic system \eqref{Linear-Hyperbolic-Syst} with the boundary condition \eqref{BC-LinearHypo}. The corresponding initial conditions are given in \eqref{IC-Linear-Hyperbolic}. By \cite{GHW-2020}, the following conclusion holds:
\begin{lemma}\label{Lemma_Local_Sol_Hyperbolic_Sys}
	Assume that
	\begin{equation*}
	  \begin{aligned}
	    \mathbf{E}_1 : = & \sum_{k=1,2} \Big\{ \|(\rho_k^{in}, u_k^{in}, \theta_k^{in})\|^2_{\mathcal{H}^{s_k} (\R^3_+)} \\
	    & + \sup_{t \in (0, \tau)} \big[\|(\mathcal{F}^\perp_\theta (f_k), \mathcal{F}^\perp_u (f_k))(t)\|^2_{\mathcal{H}^{s+1} (\R^3_+)} + \|\mathcal{J}_k (\M^0) (t)\|_{\mathcal{H}^{s+2} (\R^2)}^2 \big] \Big\} < +\infty
	  \end{aligned}
	\end{equation*}
	with $s_2 \ge 3$ and $s_1 = s_2 + 3$. Then there exists a unique smooth solution $(\rho_k, u_k, \theta_k)$ to \eqref{Linear-Hyperbolic-Syst}-\eqref{IC-Linear-Hyperbolic} with the boundary condition \eqref{BC-LinearHypo} over $t \in [0, \tau]$ such that for $k=1,2$,
	\begin{equation}\label{Energy_Inequ_LHS}
	  \begin{aligned}
	    \sup_{t \in [0, \tau]} {\|( \rho_k, u_k, \theta_k) (t)\|_{\mathcal{H}^{s_k} (\R^3_+)}^2} \le C (\tau, E_{s_1+2}) \mathbf{E}_1 \,,
	\end{aligned}
	\end{equation}
	where $E_{s_1+2} : = \sup_{t\in [0, \tau]} \| (\rho - \rho_\#, \u, T - T_\#) (t)\|_{H^{s_1+2}(\R^3_+)}$ is defined in \eqref{Compressible_Euler_Bound}.
\end{lemma}

Second, as shown in Section \ref{Sec:FormalAnaly} above, the fluid variables $(\rho_k^b, u_k^b, \theta_k^b)$ of $F_k^b$ ($k=1,2,3,4$) subject to the linear Prandtl system \eqref{Linear-Prandtl} with the boundary condition \eqref{BC-Prandtl} and the equations \eqref{u_b_k+1-derivative}-\eqref{p_b_k+1-derivative}. By the similar arguments of Lemma 2.2 in \cite{JLT-arXiv-2021} or Lemma 4.1 in \cite{GHW-2020}, one has the following result.

\begin{lemma}\label{Prop-Prandtl}
	Let $s \geq 3$, $l \geq 0$. Assume
	\begin{equation}\label{IC_norm_theta}
		\begin{aligned}
			\mathbb{E}_1 : = & \| \bar{u}^{b,in}_k \|^2_{\mathbb{H}^s_l (\R^3_+)} + \sup_{t\in [0, \tau]} \big( \| u_w (t) \|^2_{\mathbb{H}^{s+3} (\R^2)} + \| \mathtt{f}_{k-1}^b (t) \|^2_{\mathbb{H}^{s+1}_l (\R^3_+)} \big) \\
			& + \| \theta_k^{b, in} \|^2_{\mathbb{H}^s_l (\R^3_+)} + \sup_{t\in [0, \tau]} \big( \| T_w (t) \|^2_{\mathbb{H}^{s+3} (\R^2)} + \| \mathtt{g}_{k-1}^b (t) \|^2_{\mathbb{H}^{s+1}_l (\R^3_+)} \big) < \infty \,.
		\end{aligned}
	\end{equation}
	Then there exists a unique smooth solution $( \bar{u}_k, \theta_k) (t, \bar{x}, \zeta)$ to \eqref{Linear-Prandtl}-\eqref{IC-Prandtl}-\eqref{BC-Prandtl} over $t \in [0, \tau]$ satisfying
	\begin{equation}\label{theta-bnd}
		\begin{aligned}
			\sup_{t\in [0, \tau]} \Big( \| ( \bar{u}_k, \theta_k) & (t) \|^2_{\mathbb{H}^s_l (\R^3_+)} + \int_0^t \| \p_{\zeta} ( \bar{u}_k, \theta_k) (t') \|^2_{\mathbb{H}^s_l (\R^3_+)} \d t' \Big) \leq C (\tau, E_{s+1}) \mathbb{E}_1 \,,
		\end{aligned}
	\end{equation}
	where $E_{s+1}$ is defined in \eqref{Compressible_Euler_Bound}.
\end{lemma}

Third, the Knudsen boundary layers $F^{bb}_k$ will be controlled. Let
\begin{align*}
  w_{a, \beta} = \l v \r^\beta (\M^0)^{-a}
\end{align*}
with $\beta \ge 3$, $0 < a < \frac{1}{4}$. Denote by
\begin{equation*}
  \begin{aligned}
    \mathscr{S}_k : = S^{bb}_{k, 2} - \big( v_3 \p_\xi f^{bb}_{k,1} - S^{bb}_{k, 1} + \mathcal{L}^0 f^{bb}_{k,1} \big) \in (\mathcal{N}^0)^\perp \,.
  \end{aligned}
\end{equation*}
From Huang and Wang's work \cite{Huang-Wang-arXiv-2021}, the following pointwise estimates to the linear Knudsen boundary layer equation \eqref{f_bb_k_2} hold:
\begin{lemma}\label{Lemma_Linear_Higher_Knudsen}
  Assume $\int_{v \cdot n < 0} (v \cdot n) \mathbbm{f}_k (t, \bar{x}, v) \sqrt{\M^0} \d v =0$. Let $\sigma_0 \in (0,1)$ be suitably small. For any fixed integer $m \geq 0$, if
  \begin{align}\label{LKL_E_m}
  	\mathcal{E}_m := \sup_{(t, \bar{x}) \in [0, \tau] \times \R^3_+} \sum_{|\alpha| \le m} \big( \| e^{\sigma_0 \xi} \nu^{-1} w_{a, \beta} \p_{t, \bar{x}}^\alpha \mathscr{S}_k \|_{L^\infty_{\xi, v}} + \| w_{a, \beta} \p_{t, \bar{x}}^\alpha \mathbbm{f}_k \|_{L^\infty_v} \big) < \infty \,,
  \end{align}
  then the linear Knudsen boundary layer equation \eqref{f_bb_k_2} admits a unique solution $f^{bb}_{k,2} (t, \bar{x}, \xi, v)$ satisfying
  \begin{equation}\label{Linear_Knudsen_Eq_Est_Higher}
  	\sup_{(t, \bar{x}) \in [0, \tau] \times \R^3_+} \sum_{|\alpha| \le m} \big( \| e^{\sigma \xi}  w_{a, \beta} \p_{t, \bar{x}}^\alpha f^{bb}_{k,2} \|_{L^\infty_{\xi, v}} + \| w_{a, \beta} \p_{t, \bar{x}}^\alpha f^{bb}_{k,2} |_{\xi =0} \|_{L^\infty_v} \big)
  	\le \tfrac{C}{\sigma_0 - \sigma} \mathcal{E}_m
  \end{equation}
  for some constant $C > 0$ and any $\sigma \in (0, \sigma_0)$.
\end{lemma}

Recall the relations \eqref{Relt1}, \eqref{KinPart-12}, \eqref{BC-LinearHypo}, \eqref{Relt2}, \eqref{Relt3} and \eqref{Relt4}. Together with the standard Sobolev theory, the following proposition is immediately implied by Lemma \ref{Lmm-fbb-k1}, \ref{Lemma_Local_Sol_Hyperbolic_Sys} and \ref{Lemma_Linear_Higher_Knudsen}. Similar arguments can be also found in \cite{GHW-2020}.

\begin{proposition}\label{Prop-HigherOrd}
	Let $\beta \geq 3$ and $0 < \tfrac{1}{4 z} (1 - z) < a < \tfrac{1}{4}$, where $z \in (\tfrac{1}{2}, 1)$ is given in \eqref{M-Bound}. Let integers $s_0$ in Proposition \ref{Proposition_Compressible_Euler} be given large enough. $(\rho, \u, T)$ is constructed in Proposition \ref{Proposition_Compressible_Euler} and generates the local Maxwellian $\M$. There are large enough $s_k, s_k^b, s_k^{bb}, l_k^b \in \mathbb{N}_+$, $p_k, p_k^b, p_k^{bb} \in \R_+$ $(1 \leq k \leq 6)$ such that if $(\rho^{in}_k, u^{in}_k, \theta^{in}_k)$ in \eqref{IC-Linear-Hyperbolic} satisfies 
	$$\mathcal{E}_1^{in} = \sum_{1 \leq k \leq 4} \| (\rho^{in}_k, u^{in}_k, \theta^{in}_k) \|_{\mathcal{H}^{s_k} (\R^3_+)} < \infty $$ 
	and $ (\bar{u}^{b,in}_k, \theta^{b, in}_k ) $ in \eqref{IC-Prandtl} obeys $$\mathcal{E}_2^{in} = \sum_{1 \leq k \leq 4} \| (\bar{u}^{b, in}_k, \theta_k^{b, in}) \|_{\mathbb{H}^{s_k^b}_{l^b_k} (\R^3_+)} < \infty \,,$$ 
	then there are solutions $F_k = \sqrt{\M} f_k$, $F^b_k = \sqrt{\M^0} f^b_k$ and $F^{bb}_k = \sqrt{\M^0} f^{bb}_k$ $(1 \leq k \leq 6)$ constructed in Section \ref{Sec:FormalAnaly} over the time interval $t \in [0, \tau]$ with the bound
	\begin{equation}
	  \begin{aligned}
	    \sup_{t \in [0, \tau]} \sum_{k=1}^6 \Big\{ & \sum_{\gamma + |\beta| \leq s_k} \| \l v \r^{p_k} \M^{-a} \p_t^\gamma \p_x^\beta f_k (t) \|_{L^2_x L^\infty_v} \\
	    & + \sum_{j=0}^{s_k^b} \sum_{2 \gamma + |\bar{\beta}| = j} \| \l v \r^{p_k^b} (\M^0)^{-a} \p_t^\gamma \p_{\bar{x}}^{\bar{\beta}} f^b_k (t) \|_{L^2_{l^k_j} L^\infty_v} \\
	    & + \sum_{\gamma + |\bar{\beta}| \leq s_k^{bb}} \| e^{ \frac{\xi}{2^{k-1}}} \l v \r^{p_k^{bb}} (\M^0)^{-a} \p_t^\gamma \p_{\bar{x}}^{\bar{\beta}} f_k^{bb} (t) \|_{L^\infty_{\bar{x}, \xi, v} \cap L^2_{\bar{x}} L^\infty_{\xi, v}} \Big\} \\
	    \leq & C \Big(\tau, E_{s_0}, \mathcal{E}^{in}_1, \mathcal{E}^{in}_2 \Big) \,.
	  \end{aligned}
	\end{equation}
\end{proposition}

We remark that, for $1 \leq k \leq 6$, the viscous boundary layers $F_k^b$ decay algebraically associated with $\zeta > 0$, and the Knudsen boundary layers $F_k^{bb}$ decay exponentially associated with $\xi > 0$. These suffice to dominate the quantities $R_\eps^b$ and $R_\eps^{bb}$ (in \eqref{R_eps_b} and \eqref{R_eps_bb}, respectively) while deriving the uniform $L^2$-$L^\infty$ bounds for the remainder $F_{R,\eps}$.

\section{Uniform bounds for remainder $F_{R, \eps}$: Proof of Theorem \ref{Main-Thm}}\label{Sec:Remainder}

In this section, the $L^2$-$L^\infty$ arguments \cite{Guo-2010-ARMA,Guo-Jang-Jiang-2010-CPAM} will be applied to prove the main theorem. This approach is sufficient to estimate norms $\| f_{R,\eps} (t) \|_2$ and $\| h_{R,\eps}^\ell (t) \|_\infty$ for the remainder $F_{R,\eps}$. Here $f_{R,\eps}$ and $h_{R,\eps}^\ell$ are given in \eqref{f_Reps-h_Reps}, i.e.,
\begin{equation*}
  \begin{aligned}
    f_{R, \eps} = \tfrac{F_{R, \eps}}{\sqrt{\M}} \,, \quad h_{R, \eps}^\ell = \l v \r^\ell \tfrac{F_{R,\eps}}{\sqrt{\M_M}} \,.
  \end{aligned}
\end{equation*}
The proof relies on an interplay between $L^2$ and $L^\infty$ estimates for the Boltzmann equation.

Remark that the known coefficients $\M$, $F_k$, $F^b_k$, $F^{bb}_k$ ($1 \leq k \leq 6$) and the quantities $R_\eps$, $R^b_\eps$, $R^{bb}_\eps$ composed of $\M$, $F_k$, $F^b_k$, $F^{bb}_k$ are involved in the remainder equation \eqref{Remainder_Eq} with the complete diffusive boundary condition \eqref{Remainder_BC}. Thanks to Proposition \ref{Proposition_Compressible_Euler} and \ref{Prop-HigherOrd}, these knowns can be bounded by some harmless constants $C > 0$ independent of $\eps$ in the $L^2$-$L^\infty$ arguments, see also the similar arguments in Section 6 of \cite{GHW-2020}. For simplicity, the corresponding details will be omitted.

There are two main lemmas to finish the proof of Theorem \ref{Main-Thm}.

\begin{lemma}[$L^2$ Estimates]\label{Lmm-L2}
	Under the same assumptions in Proposition \ref{Prop-HigherOrd} and $\ell \geq 7$, let $(\rho, \u, T)$ be a smooth solution to the Euler equations \eqref{Compressible_Euler_Sys} over $t \in [0, \tau]$ obtained in Proposition \ref{Proposition_Compressible_Euler}. Let $c_0 > 0$ be given in \eqref{Hypocoercivity}, and the wall Maxwellian $M_w$ be assumed in Theorem \ref{Main-Thm}. Then there are constants $\eps'_0 > 0$ and $C = (\M, F_k, F^b_k, F^{bb}_k; 1 \leq k \leq 6) > 0$ such that for all $0 < \eps < \eps_0'$ and $t \in [0, \tau]$,
	\begin{equation}
	  \begin{aligned}
	    \tfrac{1}{2} \tfrac{\d}{\d t} \| f_{R,\eps} (t) \|_2^2 + \tfrac{c_0}{2 \eps} \| (\mathcal{I} - & \mathcal{P}) f_{R, \eps} (t) \|_\nu^2 +\tfrac{1}{2} \iint_{\Sigma_+} |v \cdot n| \big| (\mathcal{I} - \mathcal{D}_w) f_{R, \eps} (t) \big|^2 \d v \d \bar{x} \\
	    \le & C \sqrt{\eps} \| \sqrt{\eps}^3 h^\ell_{R, \eps} (t) \|_{\infty} \| f_{R, \eps} (t) \|_2 + C \| f_{R, \eps} (t) \|_2 \big( 1 + \|f_{R,\eps} (t) \|_2 \big)\,.
	  \end{aligned}
	\end{equation}
\end{lemma}

\begin{lemma}[$L^\infty$ Estimates]\label{Lmm-Linfty}
	Under the same assumptions in Lemma \ref{Lmm-L2}, there are constants $\eps''_0 > 0$ and $C = (\M, F_k, F^b_k, F^{bb}_k; 1 \leq k \leq 6) > 0$ such that for all $0 < \eps < \eps_0''$ and $t \in [0, \tau]$,
	\begin{equation}\label{h-infty-est}
	  \sup_{s \in [0, t]} \| \sqrt{\eps}^3 h^\ell_{R, \eps} (s) \|_{\infty} \le C \big( \| \sqrt{\eps}^3 h^\ell_{R, \eps} (0) \|_{\infty} + \sup_{0 \le s \le \tau} \| f_{R, \eps} (s) \|_2 + \sqrt{\eps}^5 \big) \,.
	\end{equation}
\end{lemma}

\begin{proof}[Proof of Theorem \ref{Main-Thm}]
	From Lemma \ref{Lmm-L2} and \ref{Lmm-Linfty},
	\begin{equation*}
	  \begin{aligned}
	    \tfrac{\d}{\d t} \| f_{R,\eps} (t) \|_2 \leq C \sqrt{\eps} ( \| \sqrt{\eps}^3 h^{\ell}_{R, \eps} (0) \|_\infty + \sup_{0 \le s \le \tau} \| f_{R, \eps} (s) \|_2 + \sqrt{\eps}^5 ) + C ( 1 + \| f_{R,\eps} \|_2 )
	  \end{aligned}
	\end{equation*}
	holds for all $0 < \eps < \min \{ \eps_0', \eps_0'' \}$ and $t \in [0, \tau]$. The Gr\"onwall inequality yields
	\begin{equation*}
	  \begin{aligned}
	    \| f_{R,\eps} (t) \|_2 + 1 \leq (\| f_{R,\eps} (0) \|_2 + 1) e^{C t} + \sqrt{\eps} e^{C t} ( \| \sqrt{\eps}^3 h_{R,\eps} (0) \|_\infty + \sup_{t \in [0, \tau]} \| f_{R,\eps} (t) \|_2 + \sqrt{\eps}^5 )  \,.
	  \end{aligned}
	\end{equation*}
	Choosing $\eps_0 : = \min \{ \eps_0', \eps_0'', \tfrac{1}{4} e^{- 2 C \tau} \} > 0$, there holds
	\begin{equation*}
	\begin{aligned}
	  \sup_{t \in [0, \tau]} \| f_{R,\eps} (t) \|_2 \leq C_\tau \big( 1 + \| f_{R, \eps} (0) \|_2 + \| \sqrt{\eps}^3 h_{R,\eps}^\ell (0) \|_\infty \big)
	\end{aligned}
	\end{equation*}
	for some $C_\tau > 0$ and any $0 < \eps < \eps_0$, which, together with \eqref{h-infty-est}, finishes the proof of Theorem \ref{Main-Thm}.
\end{proof}

\subsection{$L^2$ estimates for $f_{R,\eps}$: Proof of Lemma \ref{Lmm-L2}}

From the definition of $f_{R,\eps}$ in \eqref{f_Reps-h_Reps} and the equation of $F_{R, \eps}$ in \eqref{Remainder_Eq} with the diffusive boundary condition \eqref{Remainder_BC}, $f_{R,\eps}$ subjects to
\begin{equation}\label{Remainder_f}
 \begin{aligned}
   \p_t f_{R, \eps} + v \cdot \nabla_x & f_{R, \eps} + \tfrac{1}{\eps} \mathcal{L} f_{R, \eps}
   \\
   = & - \tfrac{(\p_t + \bar{x} \cdot \nabla_{\bar{x}}) \sqrt{\M}}{\sqrt{\M}} f_{R, \eps} + \sqrt{\eps}^3 \tfrac{1}{\sqrt{\M}} B (\sqrt{\M} f_{R, \eps}, \sqrt{\M} f_{R, \eps}) + \tfrac{1}{\sqrt{\M}} (R_\eps + R^{bb}_\eps) \\
   & + \sum_{i=1}^6 \sqrt{\eps}^{i-1} \tfrac{1}{\sqrt{\M}} \big[ B ( F_i + F^{bb}_i , \sqrt{\M} f_{R, \eps}) + B ( \sqrt{\M} f_{R, \eps}, F_i + F^{bb}_i) \big]
 \end{aligned}
\end{equation}
with
\begin{equation}\label{Remainder_BC_f}
  \gamma_- f_{R, \eps} = \mathcal{D}_w f_{R, \eps} \ \textrm{ on } \Sigma_- \,,
\end{equation}
and initial data
\begin{equation}\label{Remainder_IC_f}
 f_{R, \eps} (0, x, v) = f^{in}_{R, \eps} (x, v) : = \tfrac{F_{R, \eps}^{in} (x,v)}{\sqrt{\M^{in}}} \,.
\end{equation}
Taking inner product with $f_{R, \eps}$ with respect to $(x, v)$ on $\R^3_+ \times \R^3$, one has
\begin{equation}\label{Remainder_L2-1}
  \begin{aligned}
    \tfrac{1}{2} \tfrac{\d}{\d t} \| f_{R,\eps} \|_2^2 +  \tfrac{1}{\eps} \l \mathcal{L} f_{R, \eps}, f_{R, \eps} \r - \tfrac{1}{2} \iint_{\Sigma} v_3 \big| f_{R, \eps} (t, \bar{x}, 0, v) \big|^2 \d v \d \bar{x} \\
    = - \iint_{\R^3_+ \times R^3} \tfrac{(\p_t + v \cdot \nabla_{x}) \sqrt{\M}}{\sqrt{\M}} |f_{R, \eps}|^2 \d v \d x + \sqrt{\eps}^3 \iint_{\R^3_+ \times R^3} \tfrac{1}{\sqrt{\M}} B(\sqrt{\M} f_{R, \eps}, \sqrt{\M} f_{R, \eps}) f_{R, \eps} \d v \d x \\
    + \sum_{i=1}^6 \sqrt{\eps}^{i-1} \iint_{\R^3_+ \times R^3} \tfrac{1}{\sqrt{\M}} \big[ B (F_i + F^{bb}_i, \sqrt{\M} f_{R, \eps}) + B( \sqrt{\M} f_{R, \eps}, F_i + F^{bb}_i) \big] f_{R, \eps} \d v \d x \\
    + \iint_{\R^3_+ \times R^3} \tfrac{1}{\sqrt{\M}} (R_\eps + R^{bb}_\eps) f_{R, \eps} \d v \d x\,.
  \end{aligned}
\end{equation}

Following the arguments of \cite{Guo-Jang-Jiang-2010-CPAM}, the terms on the right-hand side of \eqref{Remainder_L2-1} can be bounded by
\begin{equation}\label{Remainder_L2_RHS}
  \begin{aligned}
    {\text{RHS\ of\ \eqref{Remainder_L2-1}}}
    \le C_\lambda \eps^2 \| h^\ell_{R, \eps} \|_{\infty} \| f_{R, \eps} \|_2 + C \| f_{R, \eps} \|_2 (1 + \| f_{R, \eps} \|_2) + \tfrac{C \lambda^2}{\eps} \| (\mathcal{I} - \mathcal{P}) f_{R, \eps} \|_{\nu}^2
  \end{aligned}
\end{equation}
for some small $\lambda > 0$ to be determined, where $C_\lambda >0$ is a constant depending on $\lambda$. Furthermore, the inequality \eqref{Hypocoercivity} reduces to
\begin{equation}\label{L2-hypercoer}
  \begin{aligned}
    \tfrac{1}{\eps} \l \mathcal{L} f_{R, \eps}, f_{R, \eps} \r \geq \tfrac{c_0}{\eps} \| (\I - \P) f_{R, \eps} \|^2_\nu \,.
  \end{aligned}
\end{equation}
It remains to deal with the boundary integral term $- \tfrac{1}{2} \iint_{\Sigma} v_3 \big| f_{R, \eps} (t, \bar{x}, 0, v) \big|^2 \d v \d \bar{x}$ in \eqref{Remainder_L2-1}. In fact, there holds
\begin{align}\label{BC-integral-f}
  - \tfrac{1}{2} \iint_{\Sigma} v_3 |f_{R, \eps} (t, \bar{x}, 0, v)|^2 \d v \d \bar{x} = \tfrac{1}{2} \iint_{\Sigma_+} |v \cdot n| \big| (\mathcal{I} - \mathcal{D}_w) f_{R, \eps} \big|^2 \d v \d \bar{x}
\end{align}
under the diffusive boundary condition \eqref{Remainder_BC_f}. Indeed, by \eqref{Remainder_BC_f} and $v \cdot n = - v_3$,
\begin{align*}
  & - \tfrac{1}{2} \iint_{\Sigma} v_3 |f_{R, \eps} (t, \bar{x}, 0, v)|^2 \d v \d \bar{x} \\
  = & \tfrac{1}{2} \iint_{\Sigma_+} v \cdot n |f_{R, \eps} (t, \bar{x}, 0, v)|^2 \d v \d \bar{x} + \tfrac{1}{2} \iint_{\Sigma_-} v \cdot n |\mathcal{D}_w f_{R, \eps} (t, \bar{x}, 0, v)|^2 \d v \d \bar{x} \,.
\end{align*}
Recalling the definition of the diffusive boundary operator $\mathcal{D}_w$ in \eqref{Dw}, the integral $\int_{v \cdot n < 0} (v \cdot n) |\mathcal{D}_w f_{R, \eps} (t, \bar{x}, 0, v)|^2 \d v $ can be calculated by
\begin{equation*}
  \begin{aligned}
    & \int_{v \cdot n < 0} (v \cdot n) |\mathcal{D}_w f_{R, \eps} (t, \bar{x}, 0, v)|^2 \d v \\
    = & \Big( \int_{v^\prime \cdot n >0} (v^\prime \cdot n) f_{R, \eps} (t, \bar{x}, 0, v^\prime) \sqrt{\M^0} \d v^\prime \Big)^2 \int_{v \cdot n < 0} (v \cdot n ) \tfrac{M_w^2}{\M^0} \d v \\
    = & - \tfrac{1}{\rho_0} \sqrt{\tfrac{2 \pi}{T^0}} \Big( \int_{v^\prime \cdot n >0} (v^\prime \cdot n) f_{R, \eps} (t, \bar{x}, 0, v^\prime) \sqrt{\M^0} \d v^\prime \Big)^2 \,,
  \end{aligned}
\end{equation*}
where one has used the fact $\int_{v \cdot n < 0} (v \cdot n ) \frac{M_w^2}{\M^0} \d v = \tfrac{1}{\rho_0} \sqrt{\tfrac{2 \pi}{T^0}}$ derived from the relations $\frac{M_w}{\M^0} = \frac{1}{\rho_0} \sqrt{\frac{2 \pi}{T^0} }$ and \eqref{Mw-Condition}, i.e., $\int_{v \cdot n < 0} |v \cdot n| M_w \d v = 1$.

Next, one focuses on the integral $\int_{v \cdot n > 0} (v \cdot n) |f_{R, \eps} (t, \bar{x}, 0, v)|^2 \d v$. By the decomposition $f_{R, \eps} = (\mathcal{I} - \mathcal{D}_w) f_{R, \eps} + \mathcal{D}_w f_{R, \eps}$,

\begin{equation*}
\begin{aligned}
  & \int_{v \cdot n >0} (v \cdot n) |f_{R, \eps} (t, \bar{x}, 0, v)|^2 \d v \\
  = & \int_{v \cdot n >0} (v \cdot n) \big| (\mathcal{I} - \mathcal{D}_w) f_{R, \eps} (t, \bar{x}, 0, v) \big|^2 \d v + \underbrace{\int_{v \cdot n >0} (v \cdot n) \big| \mathcal{D}_w f_{R, \eps} (t, \bar{x}, 0, v) \big|^2 \d v }_{B_1} \\
  & + 2 \underbrace{ \int_{v \cdot n >0} (v \cdot n) \big( \mathcal{I} - \mathcal{D}_w \big) f_{R, \eps} (t, \bar{x}, 0, v) \cdot \mathcal{D}_w f_{R, \eps} (t, \bar{x}, 0, v) \d v }_{B_2} \,.
\end{aligned}
\end{equation*}
By $\u^0 \cdot n = 0$ and $(u_w, T_w) = (\u^0, T^0)$, one has $\int_{\R^3} (v \cdot n) M_w \d v = 0$, which, combining with $\int_{v \cdot n < 0} (v \cdot n) M_w \d v = - 1$ given in \eqref{Mw-Condition}, implies that $\int_{v \cdot n > 0} (v \cdot n) M_w \d v = 1$. Together with the definition of $\mathcal{D}_w$ in \eqref{Dw} and $\frac{M_w}{\M^0} = \frac{1}{\rho_0} \sqrt{\frac{2 \pi}{T^0} }$,
\begin{align*}
  B_1 = & \int_{v \cdot n >0} (v \cdot n) \Big\{ \tfrac{M_w}{\sqrt{\M^0}} \int_{v^\prime \cdot n >0} (v^\prime \cdot n) f_{R, \eps} (t, \bar{x}, 0, v^\prime) \sqrt{\M^0} \d v^\prime \Big\}^2 \d v \\
  = & \tfrac{1}{\rho^0} \sqrt{\tfrac{2 \pi}{T^0}} \int_{v \cdot n >0} (v \cdot n) M_w \d v \Big( \int_{v^\prime \cdot n >0} (v^\prime \cdot n) f_{R, \eps} (t, \bar{x}, 0, v^\prime) \sqrt{\M^0} \d v^\prime \Big)^2 \\
  = & \tfrac{1}{\rho^0} \sqrt{\tfrac{2 \pi}{T^0}} \Big( \int_{v^\prime \cdot n >0} (v^\prime \cdot n) f_{R, \eps} (t, \bar{x}, 0, v^\prime) \sqrt{\M^0} \d v^\prime \Big)^2 \,.
\end{align*}
For $B_2$,
\begin{align*}
  B_2 = & \int_{v \cdot n >0} (v \cdot n) f_{R, \eps} (t, \bar{x}, 0, v) \mathcal{D}_w f_{R, \eps} (t, \bar{x}, 0, v) \d v - \int_{v \cdot n >0} (v \cdot n) \big| \mathcal{D}_w f_{R, \eps} (t, \bar{x}, 0, v) \big|^2 \d v \\
  := & B_{21} - B_1 \,.
\end{align*}
The goal is to show that $B_{21}=B_1$, which means $B_2 = 0$. Actually, by $\frac{M_w}{\M^0} = \frac{1}{\rho_0} \sqrt{\frac{2 \pi}{T^0} }$,
\begin{align*}
B_{21} = & \int_{v \cdot n > 0} (v \cdot n) f_{R, \eps} (t, \bar{x}, 0, v) \tfrac{1}{\rho^0} \sqrt{\tfrac{2\pi}{T^0}} \sqrt{\M^0} \int_{v^\prime \cdot n >0} (v^\prime \cdot n) f_{R, \eps} (t, \bar{x}, 0, v' ) \sqrt{\M^0} \d v^\prime \d v \\
= & \tfrac{1}{\rho^0} \sqrt{\tfrac{2\pi}{T^0}} \Big( \int_{v^\prime \cdot n >0} (v^\prime \cdot n) f_{R, \eps} (t, \bar{x}, 0, v') \sqrt{\M^0} \d v^\prime \Big)^2 = B_1 \,.
\end{align*}
Thus $B_2 = 0$. Observing $\Sigma_+ = \R^2 \times \{ v: v \cdot n > 0 \}$ and collecting the above calculations, one concludes the claim \eqref{BC-integral-f}.

Consequently, from \eqref{Remainder_L2-1}, \eqref{Remainder_L2_RHS}, \eqref{L2-hypercoer} and \eqref{BC-integral-f},
\begin{equation}\label{Remainder_L2-2}
  \begin{aligned}
   & \tfrac{1}{2} \tfrac{\d}{\d t} \| f_{R,\eps} \|_2^2 + \tfrac{c_0}{\eps} \| (\mathcal{I} - \mathcal{P}) f_{R, \eps} \|_\nu^2 + \tfrac{1}{2} \iint_{\Sigma_+} |v \cdot n| \big| (\mathcal{I} - \mathcal{D}_w) f_{R, \eps} \big|^2 \d v \d \bar{x} \\
   \le & C_\lambda \eps^2 \| h^\ell_{R, \eps} \|_{\infty} \| f_{R, \eps} \|_2 + C \| f_{R, \eps} \|_2 \big( 1 + \|f_{R,\eps}\|_2 \big) + \tfrac{C \lambda^2}{\eps} \| (\mathcal{I} - \mathcal{P}) f_{R,\eps} \|_\nu^2
  \end{aligned}
\end{equation}
for small $\lambda > 0$ to be determined. Taking small $\lambda >0$ such that $c_0 - C \lambda^2 \ge \frac{1}{2} c_0 >0$, one then finishes the proof of Lemma \ref{Lmm-L2}.

\subsection{$L^\infty$ estimates for $h^\ell_{R, \eps}$: Proof of Lemma \ref{Lmm-Linfty}}

\begin{figure}[h]
	\centering
	\begin{tikzpicture}[node distance=2cm]
	\draw[-, very thick] node[above left,scale=1]{$\p \R^3_+$}(0,0.5)--(4,0.5);
	\draw[->, thick] (2,0.4)--(2,-0.1) node[below]{$n$};
	\draw[-,thick, dashed] (0.5,2.5)--(1.5,1.5);
	\draw[->,thick] (1.5,1.5)--(2,1) node[right]{$v$};
	\draw (1.5,1.5) node[below left, scale=1]{$x$};
	\draw[->,dashed, color=red] (1.5,1.5)--(4,1.5) node[right, color=black]{time $t$};
	\end{tikzpicture}
	\caption{The case $v \cdot n \geq 0$: $t_b (x,v)$ does not exists.}\label{tb-NOTexist}
\end{figure}
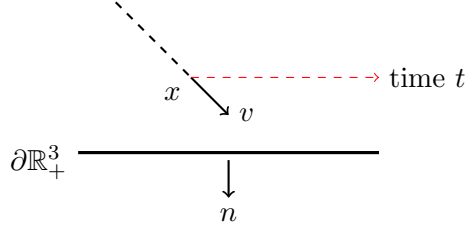

As in \cite{Guo-2010-ARMA}, the following backward characteristics of the Boltzmann equation are required to derive the $L^\infty$ estimates for $h^\ell_{R, \eps}$. Given $(t, x, v)$, define the trajectory $[X(s), V(s)]$ as
\begin{equation*}
	\left\{
	  \begin{aligned}
	   & \tfrac{\d}{\d s} X(s) = V(s)\,,\quad
	   \tfrac{\d}{\d s} V(s) = 0\,, \\
	   & [X(t), V(t)] = [x, v] \,.
	  \end{aligned}
	\right.
\end{equation*}
The solution is given by
\begin{equation*}
  [X (s), V (s)] = [X(s; t, x, v), V(s; t, x, v)] = [x- (t -s)v, v] \,.
\end{equation*}
For any fixed $(x, v) \in \overline{\R}^3_+ \times \R^3$ with $v_3 \ne 0$, its backward existed time $t_b (x , v) \ge 0$ is defined as the last moment at which the back-time straight line $[X(s; 0, x, v), V(s; 0, x, v)]$ remains in $x \in \overline{\R}^3_+$. The backward existed time can be explicitly represented by
\begin{align*}
  t_b (x, v) = \sup\big\{ \tau_0 \ge 0: x - \tau_0 v \in \R^3_+ \big\} \,,
\end{align*}
which implies $x - t_b v \in \p \R^3_+$, i.e., $x_3 - t_b v_3 =0$. Denote by
$$
x_b (x, v) = x- t_b v \in \p \R^3_+ \,.
$$
Remark here that for the case $v \cdot n <0$, the cycle will hit the boundary for one time, see Figure \ref{tb-exist}. For the case $v \cdot n \geq 0$, the back-time trajectory is a straight line and does not intersect the boundary, see Figure \ref{tb-NOTexist}.

\begin{figure}[h]
	\centering
	\begin{tikzpicture}[node distance=2cm]
	\draw[-, very thick] node[left,scale=1]{$\p \R^3_+$}(0,0)--(3,0);
    \draw[->, thick] (2,-0.1)--(2,-0.6) node[below]{$n$};
    \draw[-,thick, dashed] (1,0)--(2,1);
    \draw[->,thick] (2,1)--(2.7,1.7) node[right]{$v$};
    \draw (2,1) node[below right, scale=1]{$x$};
    \draw (1.5, 0.5) node[above left]{$t_b (x,v)$};
    \draw[->,dashed, color=red] (2,1)--(3.8,1) node[right, color=black]{time $t$};
    \draw[->,dashed, color=red] (1,0)--(3.8,0) node[right, color=black]{time $t_1 = t - t_b (x,v)$};
    \draw (1,0) node[below]{$x_b (x,v)$};
    \draw[->, thick, color=blue] (0, 0.5)--(1,0);
    \draw (0.6,0.8) node[below]{\color{blue}$v'$};
	\end{tikzpicture}
	\caption{The case $v \cdot n < 0$: $t_b (x,v)$ exists.}\label{tb-exist}
\end{figure}
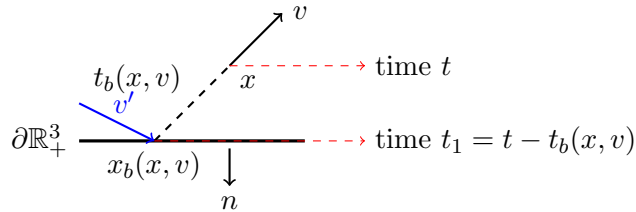

Given $x \in \overline{\R}^3_+$, $(x, v) \notin \Sigma_0 \cup \Sigma_-$ and $(t_0, v_0, x_0) = (t, x, v)$. The back-time cycle is defined as
\begin{equation}
\begin{aligned}
	X_{cl} (s; t, x, v) = & \mathbbm{1}_{[t_1, t_0)} (s) \big\{ x - (t-s) v \big\} + \mathbbm{1}_{(-\infty, t_1)} (s) \big\{ x - (t- t_1) v - (t_1 - s) v^\prime \big\} \,,\\
	V_{cl} (s; t, x, v) = & \mathbbm{1}_{[t_1, t_0)} (s) v + \mathbbm{1}_{(-\infty, t_1)} (s) v^\prime \,,
\end{aligned}
\end{equation}
where $t_1 = t - t_b (x, v)$, $v^\prime \in \mathcal{V} := \{v^\prime: v^\prime \cdot n >0\}$, see Figure \ref{tb-exist}.

From \eqref{f_Reps-h_Reps}, one has
$$
h^\ell_{R, \eps} = \tfrac{\sqrt{\M}}{\sqrt{\M_M}} \l v \r^\ell f_{R, \eps}\,.
$$
Then \eqref{Remainder_f}-\eqref{Remainder_BC_f} can be rewritten as
\begin{equation}\label{Remainder_h}
 \begin{aligned}
   \p_t h^\ell_{R, \eps} + v \cdot \nabla_x h^\ell_{R, \eps} + \tfrac{\nu (\M)}{\eps} h^\ell_{R, \eps} = q_{R,\eps}^\ell (t,x,v)
 \end{aligned}
\end{equation}
with the boundary condition
\begin{equation}\label{Remainder_BC_h}
  \gamma_- h^\ell_{R, \eps} = M_w  \tfrac{\l v \r^\ell}{\sqrt{\M_M}} \int_{v^\prime \cdot n >0} (v^\prime \cdot n) \tfrac{\sqrt{\M_M}}{\l v^\prime \r^\ell} h^\ell_{R, \eps} (t, x, v^\prime) \d v^\prime\ {\text{on}}\ \Sigma_-\,,
\end{equation}
where
\begin{equation}\label{q}
  \begin{aligned}
    q_{R,\eps}^\ell (t,x,v) = & \tfrac{1}{\eps} K_\ell h^\ell_{R, \eps} +\tfrac{\sqrt{\eps}^3 \l v \r^\ell}{\sqrt{\M_M}} B \big( \tfrac{\sqrt{\M_M}}{\l v \r^\ell} h^\ell_{R, \eps}, \tfrac{\sqrt{\M_M}}{\l v \r^\ell} h^\ell_{R, \eps} \big) + \tfrac{\l v \r^\ell}{\sqrt{\M_M}} \big( R_\eps + R^{bb}_\eps \big) \\
    & + \sum_{i=1}^6 \eps^{i-1} \tfrac{\l v \r^\ell}{\sqrt{\M_M}} \big\{ B \big( F_i + F^{bb}_i, \tfrac{\sqrt{\M_M}}{\l v \r^\ell} h^\ell_{R, \eps} \big) + B\big( \tfrac{\sqrt{\M_M}}{\l v \r^\ell} h^\ell_{R, \eps}, F_i + F^{bb}_i \big) \big\} \,,
  \end{aligned}
\end{equation}
and the operator $K_\ell g : = \l v \r^\ell K (\tfrac{g}{\l v \r^\ell})$ with
\begin{equation*}
  \begin{aligned}
    K h = & \int_{\R^3 \times \mathbb{S}^2} | (u -v) \cdot \omega | \sqrt{\M_M (u)} \tfrac{\M (v)}{\sqrt{\M_M (v)}} h (u) \d \omega \d u  \\
    & - \int_{\R^3 \times \mathbb{S}^2} | (u -v) \cdot \omega | \M (u^\prime) \tfrac{\sqrt{\M_M (v^\prime)}}{\sqrt{\M_M (v)}} h(v^\prime) \d \omega \d u \\
    & - \int_{\R^3 \times \mathbb{S}^2} | (u -v) \cdot \omega | \M (v^\prime) \tfrac{\sqrt{\M_M (u^\prime)}}{\sqrt{\M_M (v)}} h(u^\prime) \d \omega \d u \,.
  \end{aligned}
\end{equation*}

The goal of this subsection is to derive the pointwise estimates for $h_{R,\eps}^\ell (t,x,v)$ over $(t,x,v) \in [0, \tau] \times \R^3_+ \times \R^3$. According to the trajectories shown in Figure \ref{tb-exist} and \ref{tb-NOTexist}, the points set $[0, \tau] \times \R^3_+ \times \R^3$ can be clarified by the following two parts:
\begin{equation}\label{UV}
  \begin{aligned}
    \mathscr{U} = & \{ (t,x,v) \in [0, \tau] \times \R^3_+ \times \R^3 : t \in [0, \tau], v \cdot n \geq 0 \} \\
    & \qquad \cup \{ (t,x,v) \in [0, \tau] \times \R^3_+ \times \R^3 : v \cdot n < 0 \textrm{ and } t_1 = t - t_b (x,v) \leq 0 \} \,, \\
    \mathscr{V} = & [0, \tau] \times \R^3_+ \times \R^3 - \mathscr{U} \\
    = & \{ (t,x,v) \in [0, \tau] \times \R^3_+ \times \R^3 : v \cdot n < 0 \textrm{ and } t_1 = t - t_b (x,v) > 0 \} \,.
  \end{aligned}
\end{equation}
Define a number
\begin{equation}\label{t-star}
  t_\star =
  \left\{
    \begin{aligned}
      0 \,, \quad \textrm{if } (t,x,v) \in \mathscr{U} \,, \\
      t_1 \,, \quad \textrm{if } (t,x,v) \in \mathscr{V} \,.
    \end{aligned}
  \right.
\end{equation}
Observe that the trajectory over $s \in [t_\star, t]$ is
\begin{equation}\label{Trajectory}
  \begin{aligned}
    X_{cl} (s; t, x, v) = x - (t-s) v \,, \ V_{cl} (s; t, x, v) = v \,.
  \end{aligned}
\end{equation}

From integrating the equation \eqref{Remainder_h} along trajectory $[X_{cl} (s; t,x,v), V_{cl} (s; t,x,v)]$ in \eqref{Trajectory} over $s \in [t_\star , t]$,
\begin{equation}\label{h-1}
  \begin{aligned}
    h^\ell_{R, \eps} (t, x, v) = \exp\Big\{ -\tfrac{1}{\eps} \int_{t_\star}^t \nu (\phi) \d \phi \Big\} h^\ell_{R, \eps} (t_\star, X_{cl} (t_\star), V_{cl} (t_\star)) \\
    + \int_{t_\star}^t \exp\Big\{ -\tfrac{1}{\eps} \int_{t_\star}^t \nu (\phi) \d \phi \Big\} q_{R,\eps}^\ell \big( s, X_{cl} (s), V_{cl} (s) \big) \d s \,,
  \end{aligned}
\end{equation}
where $q_{R,\eps}^\ell$ is given in \eqref{q}, $t_\star$ is defined in \eqref{t-star} and the following simplified notations are employed:
\begin{equation*}
  \begin{aligned}
    [X_{cl} (s), V_{cl} (s)] : = [X_{cl} (s; t,x,v), V_{cl} (s; t,x,v)] \,, \ \nu (\phi) : = \nu (\M) \big(\phi, X_{cl} (\phi), V(\phi)\big) \,.
  \end{aligned}
\end{equation*}
Following almost the same arguments in \cite{Guo-Jang-Jiang-2010-CPAM}, the second term in the right-hand side of \eqref{h-1} can be bounded by
\begin{equation}\label{h-2}
  \begin{aligned}
    \int_{t_\star}^t \exp \Big\{ -\tfrac{1}{\eps} & \int_{t_\star}^t \nu (\phi) \d \phi \Big\} q_{R,\eps}^\ell \big( s, X_{cl} (s), V_{cl} (s) \big) \d s \\
    & \leq C (m^4 + \sqrt{\eps} + \tfrac{1}{N^0}) \sup_{s \in [0, \tau]} \| h_{R,\eps}^\ell (s) \|_\infty \\
    & + C \sqrt{\eps}^3 \sup_{s \in [0, \tau]} \| h_{R,\eps}^\ell (s) \|^2_\infty + \tfrac{C}{\sqrt{\eps}^3} \sup_{s \in [0, \tau]} \| f_{R,\eps} (s) \|_2 + C \eps
  \end{aligned}
\end{equation}
for sufficiently large $N_0 > 0$ and sufficiently small $m > 0$ to be determined. Here $C > 0$ is a harmless constant. It remains to control the first term in the right-hand side of \eqref{h-1}.

\vspace*{2mm}

{\bf Case 1:} $(t,x,v) \in \mathscr{U}$.

\vspace*{2mm}

Notice that $t_\star = 0$. Together with \eqref{nu-phi}, the first term in the right-hand side of \eqref{h-1} can be bounded by
\begin{equation}\label{h-3}
  \begin{aligned}
    \exp\Big\{ -\tfrac{1}{\eps} \int_0^t \nu (\phi) \d \phi \Big\} h^\ell_{R, \eps} (0, X_{cl} (0), V_{cl} (0)) \leq e^{ -\tfrac{C}{\eps} \l v \r t } \| h_{R,\eps}^\ell (0) \|_\infty \leq \| h_{R,\eps}^\ell (0) \|_\infty \,.
  \end{aligned}
\end{equation}
From \eqref{h-1}, \eqref{h-2} and \eqref{h-3}, one obtains that for any $(t,x,v) \in \mathscr{U}$
\begin{equation}\label{h-4}
  \begin{aligned}
    |h_{R,\eps}^\ell (t,x,v) | \leq & \| h_{R,\eps}^\ell (0) \|_\infty + C (m^4 + \sqrt{\eps} + \tfrac{1}{N^0}) \sup_{s \in [0, \tau]} \| h_{R,\eps}^\ell (s) \|_\infty \\
    & + C \sqrt{\eps}^3 \sup_{s \in [0, \tau]} \| h_{R,\eps}^\ell (s) \|^2_\infty + \tfrac{C}{\sqrt{\eps}^3} \sup_{s \in [0, \tau]} \| f_{R,\eps} (s) \|_2 + C \eps
  \end{aligned}
\end{equation}
for sufficiently large $N_0 > 0$ and sufficiently small $m > 0$ to be determined.

\vspace*{2mm}

{\bf Case 2:} $(t,x,v) \in \mathscr{V}$.

\vspace*{2mm}

Notice that $t_\star = t_1 > 0$, and $( X_{cl} (t_1), V_{cl} (t_1)) = ( x_b (x,v), v ) \in \Sigma_-$. Then the first term in the right-hand side of \eqref{h-1} is
\begin{equation*}
  \begin{aligned}
    \exp\Big\{ -\tfrac{1}{\eps} \int_{t_1}^t \nu (\phi) \d \phi \Big\} h^\ell_{R, \eps} (t_1, x_b (x,v), v) \,.
  \end{aligned}
\end{equation*}
Recalling the boundary condition \eqref{Remainder_BC_h}, one has
\begin{align*}
  h^\ell_{R, \eps} (t_1, x_b (x,v), v) = M_w \tfrac{\l v \r^\ell}{\sqrt{\M_M}} \int_{v^\prime \cdot n >0} (v^\prime \cdot n) \tfrac{\sqrt{\M_M}}{\l v^\prime \r^\ell} h^\ell_{R, \eps} (t_1, x_b (x,v), v^\prime) \d v^\prime\,.
\end{align*}
Thanks to \eqref{M-Bound} and $M_w = \tfrac{1}{\rho^0} \sqrt{\tfrac{2 \pi}{T^0}} \M^0$, one easily knows
\begin{equation*}
  \begin{aligned}
    0 < M_w \tfrac{\l v \r^\ell}{\sqrt{\M_M}} \int_{v^\prime \cdot n >0} (v^\prime \cdot n) \tfrac{\sqrt{\M_M}}{\l v^\prime \r^\ell} \d v' \leq C \,.
  \end{aligned}
\end{equation*}
Since $(t_1, x_b (x,v), v') \in \mathscr{U}$ here, $|h_{R,\eps}^\ell (t_1, x_b (x,v), v') |$ is bounded by \eqref{h-4}. Consequently, the first term in the right-hand side of \eqref{h-1} is dominated by
\begin{equation}\label{h-5}
  \begin{aligned}
    & \exp\Big\{ -\tfrac{1}{\eps} \int_{t_1}^t \nu (\phi) \d \phi \Big\} | h^\ell_{R, \eps} (t_1, x_b (x,v), v) | \\
    \leq & C \| h_{R,\eps}^\ell (0) \|_\infty + C (m^4 + \sqrt{\eps} + \tfrac{1}{N^0}) \sup_{s \in [0, \tau]} \| h_{R,\eps}^\ell (s) \|_\infty \\
    & + C \sqrt{\eps}^3 \sup_{s \in [0, \tau]} \| h_{R,\eps}^\ell (s) \|^2_\infty + \tfrac{C}{\sqrt{\eps}^3} \sup_{s \in [0, \tau]} \| f_{R,\eps} (s) \|_2 + C \eps
  \end{aligned}
\end{equation}
for sufficiently large $N_0 > 0$ and sufficiently small $m > 0$ to be determined. From \eqref{h-1}, \eqref{h-2} and \eqref{h-5}, the bound also holds for all $(t,x,v) \in \mathscr{V}$.

Finally, the proof of Lemma \ref{Lmm-Linfty} can be finished by choosing sufficiently small $\eps, m > 0$ and large $N_0 > 0$.

\appendix

\section{Vanished nonlinear Knudsen boundary layer equation}\label{Appendix}

In this appendix, we provide a proof that the leading Knudsen layer $F^{bb}_0$ must be zero. If it was supposed to exist, it formally obeys a nonlinear equation. Consider the expansion form
\begin{equation}
  \begin{aligned}
    F_\eps (t, x, v) = F_0 (t,x,v) + F^{bb}_0 (t, \bar{x}, \xi, v) + \sum_{k \geq 1} \eps^k \big[ F_k (t,x,v) + F^{bb}_k (t, \bar{x}, \xi, v) \big]
  \end{aligned}
\end{equation}
with the normal coordinate $\xi = \tfrac{x_3}{\eps}$. Plugging the previous expansion into the scaled Boltzmann equation \eqref{BE} and the complete diffusive boundary condition \eqref{Diffusive_BC}, one knows that $F^{bb}_0$ obeys
\begin{equation}
  \left\{
    \begin{aligned}
      v_3 \p_\xi F^{bb}_0 = B(F^{bb}_0, F^{bb}_0) + \big[ B(\M^0, F^{bb}_0) + B(F^{bb}_0, \M^0) \big]  \,, \\
      \gamma_- \big( \M + F^{bb}_0 \big) = M_w \int_{v^\prime \cdot n >0} (v^\prime \cdot n) \gamma_+ (\M + F^{bb}_0) \d v^\prime \,, \\
      \lim_{\xi \to +\infty} F^{bb}_0 (t, \bar{x}, \xi, v) =0\,,
    \end{aligned}
  \right.
\end{equation}
where the wall Maxwellian $M_w$ is given in \eqref{M_w2} and \eqref{uw-Tw-match}. Denote by $f^{bb}_0 = \tfrac{F^{bb}_0}{\sqrt{\M^0}}$. Then $f^{bb}_0$ subjects to
\begin{equation}\label{Knudsen_Layer_f}
  \left\{
    \begin{aligned}
      v_3 \p_\xi f^{bb}_0 + \mathcal{L}^0 f^{bb}_0 = \mathcal{Q} (f^{bb}_0, f^{bb}_0)\,, \\
      \gamma_- f^{bb}_0 = \tfrac{M_w}{\sqrt{\M^0}} \int_{v^\prime \cdot n >0} (v^\prime \cdot n) \gamma_+ f^{bb}_0 \sqrt{\M^0} \d v^\prime + \mathfrak{R} (t, \bar{x}, v)\,, \\
      \lim_{\xi \to +\infty} f^{bb}_0 (t, \bar{x}, \xi, v) =0\,,
    \end{aligned}
  \right.
\end{equation}
where
\begin{equation*}
\begin{aligned}
  \mathcal{Q} (f^{bb}_0, f^{bb}_0) = \tfrac{1}{\sqrt{\M^0}} \big[ B(\sqrt{\M^0} f^{bb}_0, \M^0) + B(\M^0, \sqrt{\M^0} f^{bb}_0) \big]
\end{aligned}
\end{equation*}
and
\begin{equation}\label{R}
\mathfrak{R} (t, \bar{x}, v) = \tfrac{M_w}{\sqrt{\M^0}} \int_{v^\prime \cdot n >0} (v^\prime \cdot n) \M^0 \d v^\prime - \sqrt{\M^0} \,.
\end{equation}
From \cite{Coron-Golse-Sulem-1988-CPAM}, the solvability condition of \eqref{Knudsen_Layer_f} is
\begin{align}\label{Solvability_Consition}
\int_{v_3 >0} v_3 \mathfrak{R} (t, \bar{x}, v) \sqrt{\M^0} \d v =0\,.
\end{align}
As in Subsection \ref{Subsec:2.2} before, the equality \eqref{Solvability_Consition} implies $\u^0 \cdot n = 0$.

Conversely, if $\u^0 \cdot n = 0$, one has $\int_{\R^3} (v \cdot n) \M^0 \d v = 0$, which means that $\int_{v \cdot n > 0} (v \cdot n) M_w \d v = - \int_{v \cdot n < 0} (v \cdot n) M_w \d v = 1$ by \eqref{Mw-Condition} and $M_w = \tfrac{1}{\rho^0} \sqrt{\tfrac{2 \pi}{T^0}} \M^0$. Consequently, there holds
\begin{equation}
  \begin{aligned}
    \mathfrak{R} (t, \bar{x}, v) = \tfrac{M_w}{\sqrt{\M^0}} \int_{v^\prime \cdot n >0} (v^\prime \cdot n) \M^0 \d v^\prime - \sqrt{\M^0} \\
    = \tfrac{1}{\rho^0} \sqrt{\tfrac{2 \pi}{T^0}} \sqrt{\M^0} \int_{v' \cdot n > 0} (v' \cdot n) \M^0 \d v - \sqrt{\M^0} \\
    = \sqrt{\M^0} \int_{v' \cdot n > 0} (v' \cdot n) M_w \d v - \sqrt{\M^0} = 0 \,,
  \end{aligned}
\end{equation}
which means that $f^{bb}_0 \equiv 0$ is a solution to \eqref{Knudsen_Layer_f}. Then the uniqueness implies that the nonlinear Knudsen boundary layer equation \eqref{Knudsen_Layer_f} only admits zero solution.


\section*{Acknowledgments}
The authors would like to thank Dr. Yulong Wu for many valuable discussions on this paper. The authors also appreciate the anonymous referees for careful reading and helpful comments. The author N. JIANG was supported by the grants from the National Natural Foundation of China under contract Nos. 11971360 and 11731008, and also supported by the Strategic Priority Research Program of Chinese Academy of Sciences, Grant No. XDA25010404. The author Y.-L. LUO is supported by grants from the National Natural Science Foundation of China under contract No. 12201220 and the Guang Dong Basic and Applied Basic Research Foundation under contract No. 2021A1515110210. The author S. J. TANG was supported by the National Natural Science Foundation of China, Grant No. 12201480 and the Fundamental Research Funds for the Central Universities grant WUT:  223114007. 
\bigskip

\bibliography{reference}

\end{document}